\documentclass[12pt,a4paper]{amsart}
\usepackage{amsmath,amssymb,amsthm, amscd}

\newtheorem{theorem}{Theorem}[section]

\newtheorem{lemma}[theorem]{Lemma}

\newtheorem{proposition}[theorem]{Proposition}

\newtheorem{corollary}[theorem]{Corollary}

\newcommand{\Fiel}{\mathbb{F}}
\newcommand{\Fq}{\mathbb{F}_{q}}
\newcommand{\Fqd}{\mathbb{F}_{q^d}}
\newcommand{\Fqbar}{\overline{\Fq}}
\newcommand{\F}{\mathbb{F}}

  \makeatletter
  \def\imod#1{\allowbreak\mkern10mu({\operator@font mod}\,\,#1)}
  \makeatother

\begin{document}
\pagestyle{plain} \pagenumbering{arabic}
\title{Real and strongly real classes in finite linear groups}

\author{Nick Gill}
\email{nickgill@cantab.net}
\address{Department of Mathematics, University Walk, Bristol, BS8 1TW, United Kingdom}

\author{Anupam Singh}
\email{anupamk18@gmail.com}
\address{IISER, Central Tower, Sai Trinity Building, near Garware Circle, Pashan, Pune 411021 India}

\thanks{Part of this work was completed when both authors were postdoctoral fellows at IMSc, Chennai. The first author is grateful to Ian Short, and others at the National University of Ireland, Maynooth, who took an interest in this work. Both authors would like to thank Rod Gow, Amritanshu Prasad and the referee for helpful comments.}

\subjclass[2000]{20E45, 20G40, 15A33}

\keywords{real, strongly real, conjugacy, general linear, special linear, finite field}

\date{\today} 
\begin{abstract}
We classify the real and strongly real conjugacy classes in $GL_n(q)$, $SL_n(q)$, $PGL_n(q)$, $PSL_n(q)$, and all quasi-simple covers of $PSL_n(q)$. In each case we give a formula for the number of real, and the number of strongly real, conjugacy classes.
\end{abstract}

\maketitle

\begin{center}
 \begin{tabular}{|ll|}
  \hline
\ref{s: intro} & Introduction \\
\ref{s: poly} & Self-reciprocal and $\zeta$-self-reciprocal polynomials \\
\ref{s: gl} & Background information and $GL_n(q)$ \\
\ref{s: sl} & $SL_n(q)$,  $n\not\equiv 2 \imod 4$ or $ q\not\equiv 3\imod 4$ \\
\ref{s: sl2} & $SL_n(q)$,  $n\equiv 2 \imod 4$  and  $q\equiv 3\imod 4$ \\
\ref{s: slstrong} & Strongly real conjugacy classes in $SL_n(q)$ \\
\ref{s: pgl} & $PGL_n(q)$ \\
\ref{s: psl1} & $PSL_n(q)$, $q$ is even or $|n|_2\neq |q-1|_2$ \\
\ref{s: psl2} & $PSL_n(q)$, $q$ is odd and $|n|_2=|q-1|_2$ \\
\ref{pslstrongly} & Strongly real classes in $PSL_n(q)$ \\
\ref{s: isogenous} & Quotients of $SL_n(q)$ \\
\ref{s: exceptional} & Some exceptional cases \\
\ref{s: small} & Some small rank calculations \\
\ref{s: further} & Further work \\
\hline
 \end{tabular}

\end{center}

\section{Introduction}\label{s: intro}

Let $G$ be a group. An element, $g$, of $G$ is called {\it real} if there exists $h\in G$ such that $hgh^{-1}=g^{-1}$. If $h$ can be chosen to be an involution (i.e. $h^2=1$) then we say that $g$ is {\it strongly real}. In all cases we say that $h$ is a {\it reversing element} for $g$. If $g$ is real (resp. strongly real) then all conjugates of $g$ are real (resp. strongly real), hence we talk about {\it real classes} and {\it strongly real classes} in $G$.

Tiep and Zalesski\cite{tz} have listed all quasi-simple groups of Lie type for which every element is real. In this paper we generalise one part of the work of Tiep and Zalesski by identifying exactly which conjugacy classes are (strongly) real in the quasi-simple groups which cover $PSL_n(q)$; furthermore we count these classes.

This paper is structured as follows: in section \ref{s: poly} we outline results concerning a special class of polynomials. In Section \ref{s: gl} we introduce some background theory that will be important for the rest of the paper. In particular we use information from Section \ref{s: poly} to classify the real and strongly real classes in $GL_n(q)$. While this work is not new, it forms the foundation for the rest of the paper. In Sections \ref{s: sl} to \ref{s: slstrong} we classify the real and strongly real classes in $SL_n(q)$.

From Section \ref{s: pgl} onwards, we focus on $PGL_n(q)$, $PSL_n(q)$ and, finally, those quasi-simple groups which cover $PSL_n(q)$. This analysis has a slightly different flavour because the groups of interest are no longer subgroups of $GL_n(q)$, but quotients of subgroups. To understand reality in $PGL_n(q)$ and $PSL_n(q)$ we need to study the {\it $\zeta$-real elements} in $GL_n(q)$ and $SL_n(q)$; the $\zeta$-real elements are defined in Section \ref{s: poly} and are studied in parallel with real elements throughout Sections \ref{s: gl} to \ref{s: sl2}. Their significance is explained by Lemma \ref{l: projreal}.

An understanding of reality in $PGL_n(q), PSL_n(q)$ and the remaining quasi-simple covers of $PSL_n(q)$ requires that we understand how conjugacy is affected when we factor out the centre of a group. This is discussed in the first half of Section \ref{s: pgl}; that discussion sets the scene for what follows in Sections \ref{s: pgl} to \ref{s: isogenous}. Sections \ref{s: gl} to \ref{s: isogenous} all include a theorem near the end which summarises the main results of the section. 

Section \ref{s: exceptional} covers some exceptional quasi-simple covers of $PSL_n(q)$ that require different techniques, and thereby completes our analysis of real and strongly real classes in the quasi-simple covers of $PSL_n(q)$. An interesting consequence of this analysis is the following statement, the proof of which is scattered throughout the paper.

\begin{theorem}
Let $G$ be isomorphic to $GL_n(q)$, $PGL_n(q)$, or a cover of $PSL_n(q)$. Then all real elements are strongly real if and only if $G$ is in the following list:
\begin{enumerate}
\item $GL_n(q)$;
\item $PGL_n(q)$;
\item $SL_n(q)/Y$ with $n\not\equiv 2 \pmod 4$ or $q$ even, here $Y$ is any central subgroup in $SL_n(q)$;
\item $SL_n(q)/Y$ with $n\equiv 2 \pmod 4$ and $q\equiv 1\pmod 4$, here $Y$ is any even order central subgroup in $SL_n(q)$;
\item $PSL_2(q)$;
\item $3.PSL_2(9)$.
\end{enumerate}
\end{theorem}

In Section \ref{s: small} we give formulae, as polynomials in $q$, for the number of real and strongly real classes in all relevant groups with $n\leq 6$. We conclude with Section \ref{s: further} in which we outline possible areas of future research.

As we have mentioned, the work on $GL_n(q)$ is not new; Gow has already enumerated the real classes for $GL_n(q)$ and given a generating function for this count \cite{gow}. The work on $SL_n(q)$ is partially new; a version of Proposition \ref{p: not2mod4}, which deals with the case when $n\not\equiv 2 \pmod 4$ or $q\not\equiv 3 \pmod 4,$ was first proved by Wonenburger \cite{wonen}. Results concerning the remaining case, when $n\equiv 2 \pmod 4$ and $q \equiv 3 \pmod 4,$ are new. As far as we know the results obtained for $PGL_n(q)$ and $PSL_n(q)$ do not exist in the literature. However Gow has communicated with us concerning work on real classes in $PGL_n(q)$; so, although this has not been published, some of the results are already known.

\section{Self-reciprocal and $\zeta$-self-reciprocal polynomials}\label{s: poly}

As we shall see, real elements in $GL_n(q)$ will turn out to correspond to sequences of self-reciprocal polynomials. In this section we define what a self-reciprocal polynomial is and we gather together some basic facts about such polynomials.

We also introduce the notion of a $\zeta$-real element in $H$, a subgroup of $GL_n(q),$ as follows: fix $\zeta$, a non-square in $\Fq$. We say that $g$ is {\it $\zeta$-real} in $H$ if $tgt^{-1}=g^{-1}(\zeta I)$ for some $t\in H$; we say that $g$ is {\it strongly $\zeta$-real} if $t$ can be taken to be an involution. Once again we say that $t$ is a {\it reversing element} for $g$. The $\zeta$-real elements of $GL_n(q)$ will turn out to be of vital importance when we come to examine the real elements of $PGL_n(q)$ and $PSL_n(q)$. Note that we will sometimes abuse notation and, for an element $g\in GL_n(q)$, write $\zeta g$ when we mean $(\zeta I)g$.

It turns out that $\zeta$-real elements will correspond to sequences of $\zeta$-self-reciprocal polynomials. Thus in this section we also examine these polynomials. Throughout what follows $\zeta$ is a fixed non-square of $\Fq$.

\subsection{Definitions}

Consider a polynomial $f(t)\in\Fq[t]$ of degree $d$ with roots $[\alpha_1,\dots,\alpha_d]$ in $\Fqbar$, the algebraic closure of $\Fq$. We say that $f(t)$ is {\it self-reciprocal} if,
$$[\alpha_1,\dots,\alpha_d]=[\alpha^{-1}_1,\dots,\alpha^{-1}_d].$$
We say that $f(t)$ is {\it $\zeta$-self-reciprocal} if
$$[\alpha_1,\dots,\alpha_d]=[\zeta\alpha^{-1}_1,\dots,\zeta\alpha^{-1}_d].$$
For both definitions, by $[,\dots, ]$ we mean an unordered list of roots, taken with multiplicity. Note that, since $\zeta$ is a non-square in $\Fq$, when we talk about $\zeta$-self-reciprocal polynomials we assume that $q$ is odd.

\subsection{Self-reciprocal polynomials}

We are interested in $T_d$, the set of self-reciprocal degree $d$ polynomials in $\Fq[t]$ with constant term equal to $1$.

It is easy enough to prove that $T_d=F_d\cup G_d$ where
$$F_d=\{f(t)=t^d+a_1t^{d-1}+a_2t^{d-2}+\dots + a_2t^2+a_1t+1\},$$
$$G_d=\{g(t)=-t^d+a_1t^{d-1}-a_2t^{d-2}+\dots + a_2t^2-a_1t+1\},$$
and the $a_i$ vary over $\Fq$. Note that if $q$ is even then $F_d$ and $G_d$ coincide. We define $n_{q,d}$ to be the number of self-reciprocal polynomials $f$ in $\mathbb F_q[t]$ of degree $d$ which satisfy $f(0)=1$.

\begin{lemma}\label{l: nqd}
The  number $n_{q,d}$ is given in the following table:
\begin{center}
\begin{tabular}{|c||c|c|}
\hline
& $q$ is odd & $q$ is even \\
\hline
\hline
$d$ is odd & $2q^{\frac{d-1}{2}}$ & $q^{\frac{d-1}2}$ \\
\hline
$d$ is even & $(q+1)q^{\frac{d}2-1}$ & $q^{\frac{d}2}$ \\
\hline
\end{tabular}
\end{center}
\end{lemma}

Before we leave self-reciprocal polynomials we make one more definition: let $p(t)=t^n+a_1t^{n-1}+\cdots+a_n$ be a monic polynomial over $k$. We define $\tilde p(t)=a_n^{-1}t^np(\frac{1}{t})$; this is the monic polynomial whose roots are precisely the inverse of the roots of $p(t)$. Thus a monic polynomial $p(t)$ is self reciprocal if and only if $p(t)=\tilde p(t)$.  Note too that $p(t)$ is irreducible in $k[t]$ if and only if $\tilde p(t)$ is irreducible in $k[t]$.

\subsection{$\zeta$-self-reciprocal polynomials}

Let $q$ be odd and let $f$ be a $\zeta$-self-reciprocal polynomial with roots in $[\alpha_1,\dots,\alpha_d]\in\Fqbar$. Suppose that
$$\alpha_i=\zeta\alpha_j^{-1}, \quad \alpha_j=\zeta\alpha_k^{-1}.$$
Then, clearly, $\alpha_i=\alpha_k$. Thus the roots of $f$ can be partitioned into subclasses of size at most $2$. Observe that, within these subclasses, $\alpha_i\alpha_j=\zeta$. Now if the subclass is of size $1$ then $\alpha_i^2=\zeta$ and so does not lie in $\Fq$. We conclude that $d$ must be even.

Now, for $d$ even, define the set $S_d$ to be the union of the following two sets:
$$\left\{f(t)=\frac1{\zeta^{\frac{d}2}}t^d+a_1\frac1{\zeta^{\frac{d}2-1}}t^{d-1}+a_2\frac1{\zeta^{\frac{d}2-2}}t^{d-2}+\dots + a_2t^2+a_1t+1\right\},$$
$$\left\{g(t)=-\frac1{\zeta^{\frac{d}2}}t^d+a_1\frac1{\zeta^{\frac{d}2-1}}t^{d-1}-a_2\frac1{\zeta^{\frac{d}2-2}}t^{d-2}+\dots + a_2t^2-a_1t+1\right\},$$
where the $a_i$ vary over $\Fq$. For $d$ odd, define $S_d$ to be empty.

\begin{lemma}\label{l: nqdzeta}
The set $S_d$ is precisely the set of $\zeta$-self-reciprocal polynomials of degree $d$ with constant term equal to $1$. Moreover, $|S_d|=n_{q,d}\sigma_d$ where $\sigma_d$ equals $1$ if $d$ is even and  $0$ otherwise.
\end{lemma}
\begin{proof}
Let $h(t)$ be a polynomial of degree $d$ and write the list of roots for $h(t)$ as $[\alpha_1,\dots, \alpha_d]$. Then, clearly, $t^d h(\frac{\zeta}{t})$ is a polynomial of degree $d$ and its list of roots is $[\zeta\alpha_1^{-1}, \dots, \zeta\alpha_d^{-1}]$. Thus $h(t)$ will be $\zeta$-self-reciprocal if and only if $h(t)$ is equal to a scalar multiple of $t^d h(\frac{\zeta}{t})$.

Examining $f(t)$ and $g(t)$ given in the form above in $S_d$ we observe that $f(t)= \frac{t^d}{\zeta^{\frac{d}{2}}}f(\frac{\zeta}{t})$, while $g(t)= \frac{t^d}{-\zeta^{\frac{d}{2}}}g(\frac{\zeta}{t})$; hence all elements of $S_d$ are indeed $\zeta$-self-reciprocal.

We must now show that $S_d$ contains all $\zeta$-self-reciprocal polynomials. Let $h(t)$ be a $\zeta$-self-reciprocal polynomial and consider the roots $[\alpha_1, \dots, \alpha_d]$ of $h(t)$ split into subclasses of size at most $2$ as described above. Now recall that the subclasses of size $1$ have form $\{\alpha\}$ where $\alpha^2=\zeta$; there are an even number of these so can join them together in pairs to ensure that the list of roots is split into subclasses, $\{\alpha_i, \alpha_j\}$ of size $2$. These subclasses either satisfy $\alpha_i\alpha_j=\zeta$ or $\alpha_i\alpha_j=-\zeta$; in the latter case $\alpha_i=-\alpha_j$.

Consider these two cases. Firstly if $\alpha_i\alpha_j=\zeta$ then we can multiply the corresponding linear factors, $t-\alpha_i$ and $t-\alpha_j$, to obtain an $\Fq$-scalar multiple of
$$\frac1\zeta t^2+at+1$$
where $a\in\Fqbar$. Alternatively, if $\alpha_i=-\alpha_j^{-1}$ and $\alpha_i\alpha_j=-\zeta$, then multiplying the corresponding linear factors yields an $\Fq$-scalar multiple of
$$-\frac1\zeta t^2+1.$$
If we multiply such pairs together we generate polynomials of the following forms:
$$f(t)=\frac1{\zeta^{\frac{d}2}}t^d+a_1\frac1{\zeta^{\frac{d}2-1}}t^{d-1}+a_2\frac1{\zeta^{\frac{d}2-2}}t^{d-2}+\dots + a_2t^2+a_1t+1,$$
$$g(t)=-\frac1{\zeta^{\frac{d}2}}t^d+a_1\frac1{\zeta^{\frac{d}2-1}}t^{d-1}-a_2\frac1{\zeta^{\frac{d}2-2}}t^{d-2}+\dots + a_2t^2-a_1t+1,$$
for some $a_i\in \Fqbar$. Now $h(t)$ is of this form and lies in $\Fqbar$; in other words, for $h(t)$, the coefficients $a_i$ lie in $\Fq$ and we have the required form.

The formula for the size of $S_d$ is an easy consequence of its definition.
\end{proof}

Before we leave $\zeta$-self-reciprocal polynomials, we make one more definition: for $p(t)$ a monic polynomial of degree $d$ in $k[t]$, define $\breve p(t)$ to be the monic polynomial which is a scalar multiple of $t^d p(\frac{\zeta}{t}).$ Clearly $p(t)$ will be $\zeta$-self-reciprocal if and only if $p(t)=\breve p(t)$.

\section{Background information and $GL_n(q)$}\label{s: gl}

We start by collecting some basic facts which we will need in the sequel. Recall that, for $g$ an element of a group $G$, we denote the centralizer of $g$ by $C_G(g)$. We define the {\it reversing group} of $g$,
$$R_G(g) = \{h\in G \mid hgh^{-1} = g \textrm{ or } hgh^{-1}=g^{-1}\}.$$
When $G\leq GL_n(\Fiel)$, for some field $\Fiel$, we define a related group: fix $\zeta$ to be a non-square in $k$ and let
$$R_{G,\zeta}(g) = \{h\in G \mid hgh^{-1} = g \textrm{ or } hgh^{-1}= \zeta g^{-1}\}.$$
It is easy to check that $R_G(g)$ and $R_{G,\zeta}(g)$ are indeed groups and that they contain $C_G(g)$, the centralizer of $g$ in $G$, as a subgroup of index at most $2$. In fact, provided $g^2\neq 1$, the index of $C_G(g)$ in $R_G(g)$ (resp. $R_{\zeta, G}(g)$) is $2$ if and only if $g$ is real (resp. $\zeta$-real) in $G$; in this case $R_G(g)\backslash C_G(g)$ (resp. $R_{G,\zeta}(g)\backslash C_G(g)$) is the set of all reversing elements for $g$.

When $G$ is a subgroup of $GL_n(\Fiel)$ we can make use of the {\it Jordan decomposition} of elements of $GL_n(\Fiel)$. We need only a few basic facts about this decomposition (more details can be found in \cite[p.20]{springer2}). Any element $g\in GL_n(\Fiel)$ can be written uniquely as $g=g_sg_u$ where $g_s$ is a semi-simple element and $g_u$ is a unipotent element. We have the following result which is a generalization of \cite[Lemma 2.2.1]{singhthakur2}:

\begin{lemma}\label{l: centralizer}
Let $g=g_sg_u$ be the Jordan decomposition of $g$ in $GL_n(\Fiel)$. Let $G$ be a subgroup of $GL_n(\Fiel)$ which contains $g$. Then $g$ is real (resp. $\zeta$-real) in $G$ if and only if $g_s$ is real (resp. $\zeta$-real) in $G$ and $g_u^{-1}$ is conjugate to $xg_ux^{-1}$ in $C_G(g_s)$ where $xg_sx^{-1}=g_s^{-1}$ (resp. $xg_s x^{-1} = \zeta g_s^{-1}$).
\end{lemma}
\begin{proof}
We prove this statement for the situation when $g$ is $\zeta$-real; the case where $g$ is real has the same proof if we simply remove all instances of $\zeta$.

Suppose that $g$ is $\zeta$-real. Then
\begin{equation}\label{e: one}
hgh^{-1}=hg_s g_uh^{-1} = (hg_sh^{-1})(hg_u h^{-1}) = \zeta g^{-1}
\end{equation}
 for some $h\in G$. Now the Jordan decomposition of $\zeta g^{-1}$ gives $(\zeta g^{-1})_s = \zeta g_s^{-1}$ and $(\zeta g^{-1})_u = g_u^{-1}$. But, since $hg_sh^{-1}$ is semi-simple and $hg_uh^{-1}$ is unipotent and they commute, we have already given a Jordan decomposition of $\zeta g^{-1}$ in (\ref{e: one}). Since this decomposition is unique we must have
$$hg_sh^{-1}= \zeta g_s^{-1} \textrm{ and } hg_uh^{-1}=g_u^{-1}.$$
This implies that $g_s$ is $\zeta$-real, and that $g_u^{-1}$ is conjugate to $hg_u h^{-1}$ in $C_G(g_s)$ (in fact, in this case, the two are equal).

Now for the converse: suppose that $h\in C_G(g_s)$ satisfies $hg_u^{-1} h^{-1} = xg_u x^{-1}$. Then
$$(h^{-1}x) g (h^{-1} x)^{-1} = h^{-1}xgx^{-1} h = h^{-1}xg_sx^{-1}x g_u x^{-1} h = \zeta g_s^{-1} g_u^{-1} = \zeta g^{-1}.$$
\end{proof}

In the next subsection we will discuss the real conjugacy classes in $GL_n(\Fiel)$, for which we will need an understanding of Jordan canonical forms. Before we embark on this discussion we mention one last fact which is closely related to the Jordan canonical forms, but is of a slightly different nature. Consider an element $g\in GL_n(q)$ where $q$ is a prime power. Write $g=g_sg_u$, the Jordan decomposition of $g$. We will be interested in applying Lemma \ref{l: centralizer} and so we will need an understanding of the structure of $C_G(g_s)$ for $G=GL_n(q)$. We describe the structure of $C_G(g_s)$ in a particular case.

\begin{lemma}
Take $g\in G=GL_n(q)$ and suppose that the characteristic polynomial of $g$ and the minimal polynomial of $g$ coincide and are equal to $f(t)^r$ where $f(t)$ is an irreducible polynomial of degree $d$ and $n=dr$. Then
$$GL_r(q^d)\cong C_G(g) \leq R_G(g)\leq GL_r(q^d).\langle \sigma\rangle$$ 
where $\sigma$ is a field automorphism of order $d$.
\end{lemma}

\begin{proof}
Let $V_{q}$ be the vector space of dimension $n$ over $\Fq$ on which $GL_n(q)$ acts naturally, and let $V_{q^d}$ be an $r$-dimensional vector space over $\Fqd$. There is a natural $\Fq$-vector space isomomorphism $\phi:V_q\to V_{q^d}$ which induces an embedding of $\Fqd$ into $G$ (as the centre of $GL(V_{q^d})\leq GL(V_q)$) such that $g$ lies in $\Fqd$.

Now suppose that $h\in GL(V_q)$ centralizes $g$, i.e. $hgh^{-1} = g$. We demonstrate that $h$ lies in $GL(V_{q^d})$. It is clear, first of all, that $h(v_1+v_2)=h(v_1)+h(v_2)$ for all $v_1,v_2\in V_{q^d}$; this follows from the linearity of the action of $h$ on $V_q$. We need to demonstrate that $h$ preserves scalar multiplication, for scalars in $\Fqd$.

Observe that $\langle g\rangle$ has a well-defined action on $V_1$ where $V_1$ is any 1-dimensional subspace of $V_{q^d}$. Again we can think of $V_1$ as a vector space over $\Fq$ or $\Fqd$; the element $g|_{V_1}$ acts as an $\Fq$-vector space endomorphism with minimal polynomial $f(t)$, and as an element of $\Fqd$ by multiplication. 

Suppose that $g|_{V_1}$ fixes a proper subspace $W$ of $V_1$. Let $v$ be a non-zero vector in $W$ and consider the elements $v, gv, g^2v, \dots$ Since these all lie in $W$ and $\dim W = c<d$, we know that there exist scalars $a_i\in \Fq$ such that
$$0=a_0 v+ a_1gv +a_2g^2v + \dots a_c g^c v = (a_0+a_1g+\dots +a_cg^c)v.$$
Since $v\neq 0$ this implies that $a_0+a_1g+\dots+a_cg^c = 0$ which contradicts the fact that $f(t)$ is the minimal polynomial of $g|_{V_1}$. Thus, for any $v\in V_1$, there is an $\Fq$-basis for $V_1$ of form $\{v, gv, g^2v, \dots, g^{d-1}v\}$.

Now take $v\in W$ and $\alpha\in\Fqd$. Note that $\alpha$ commutes with $g$ since they both lie in $\Fqd$. Then
$$\alpha v = (b_0v +b_1 gv + b_2g^2v+ \dots b_{d-1}g^{d-1}v),$$
for some $b_0, \dots, b_{n-1}\in\Fq$. This implies that, for $i=1,\dots, d-1,$
$$\alpha g^iv = g^i \alpha v = g^i(b_0v +b_1 gv + b_2g^2v+ \dots b_{d-1}g^{d-1}v)
=(b_0I+b_1g + b_2g^2+\dots b_{d-1}g^{d-1})g^iv.$$
Thus $\alpha = b_0I+b_1g + b_2g^2+\dots b_{n-1}g^{d-1}.$ Now observe that, for $h\in C_G(g)$, $v\in V$,
$$h(\alpha v) = h(b_0I+b_1g + b_2g^2+\dots b_{d-1}g^{d-1})v = (b_0I+b_1g + b_2g^2+\dots b_{d-1}g^{d-1})h(v) = \alpha h(v).$$

We conclude that $C_G(g)\leq GL(V_{q^d})$. It follows immediately that $C_G(g)=GL(V_{q^d})\cong GL_r(q^d)$. 

Now we wish to study the normalizer of $C_G(g)$. Write $d=d_1\cdots d_l$ where $d_1,\dots, d_l$ are primes. Then $GL(V_{q^d}).\langle \delta\rangle$ is a maximal subgroup of $GL(V_{q^{d/d_1}})\langle \delta_1\rangle$ where $\delta$ (resp. $\delta_1$) is a field automorphism of $GL(V_{q^d})$ (resp. $GL(V_{q^{d/d_1}})$) of order $d$ (resp. $d/d_1$). (Details can be found in \cite[\S 4.3]{kl}; see, in particular, p.116.) Similarly $GL(V_{q^{d/d_1}})\langle \delta_1\rangle$ is a maximal subgroup of $GL(V_{q^{d/(d_1d_2)}})\langle \delta_2\rangle$ where $\delta_2$ is a field automorphism of $GL(V_{q^d})$ of order $d/(d_1d_2)$, and so on. We conclude that  $C_G(g)$ is normal in $GL_r(q^d).\langle \sigma \rangle$ where $\sigma$ is a field automorphism of $GL(W)$ of order $d$.

Observe that $N_{\Gamma L_n(q)}(C_G(g))$ must, therefore, contain $GL_r(q^d).\langle \overline{\sigma} \rangle$ where $\overline{\sigma}$ is a field automorphism of $GL_r(q^d)$ of order $d\log_pq$. Now any element of $\Gamma L_r(q^d)$ which normalizes $C_G(g)$ must normalize $\Fqd=Z(C_G(g))$ and so must induce a field automorphism on $\Fqd$; these are all accounted for and so we conclude that $N_{\Gamma L_n(q)}(C_G(g))=GL_r(q^d).\langle \overline{\sigma} \rangle$. It follows, therefore, that $N_{G}(C_G(g))=GL_r(q^d).\langle \sigma \rangle$ as required.
\end{proof}

Note that, in the language of Jordan canonical forms, $g$ is semi-simple and conjugate to a Jordan block matrix. Note too that we could replace $R_G(g)$ with $R_{G,\zeta}(g)$ in the statement of the lemma and it would remain true. Finally observe that, for $g$ to be real (resp. $\zeta$-real) in $GL_n(q),$ we must have $d$ even and $x(g)=g^{-1}$ (resp. $x(g) = (\zeta I) g^{-1}$) where $x$ is a field automorphism of $GL_r(q^d)$ of order $2$.

\subsection{Real Conjugacy Classes in $GL_n(\Fiel)$}\label{s: glnk}

Let $\Fiel$ be a field. The conjugacy classes in $GL_n(\Fiel)$ are determined using the theory of Jordan canonical forms. We will assume a basic understanding of this theory in what follows (more details can be found in \cite{jacobson2}). The theory is based upon the idea that, given an element $g$ of $GL_n(\Fiel)$, we can create an $n$-dimensional $\Fiel[t]$ module $V$ by defining a scalar multiplication,
$$t.v = gv, v\in V,$$
and extending linearly. The isomorphism classes of $V$ constructed in this way are in one-to-one correspondence with the conjugacy classes of $GL_n(\Fiel)$. They are also in one-to-one correspondence with the set of all multi-sets of form
$$\{f_1(t)^{a_1}, \dots, f_r(t)^{a_r}\}$$
where, for $i=1,\dots, r$, $f_i(t)$ is a monic irreducible polynomial in $\Fiel[t]$ which is not equal to $t$, $a_i$ is a positive integer, and $\sum_{i=1}^r\textrm{deg}(f_i)a_i=n$. These correspondences allow us to classify conjugacy in $GL_n(q)$.

Before we state the main result that we shall need, we introduce some notation: a {\it partition} $\nu$ of $n$ is a finite multi-set of positive integers, $\nu=\{\nu_1,\dots, \nu_r\}$ that sums to $n$; we say $|\nu|=n$). Write $\nu = 1^{n_1}2^{n_2}3^{n_3}\cdots$ to mean that
$$n=\underbrace{1+\cdots+1}_{n_1}+\underbrace{2+\cdots+2}_{n_2}+\cdots$$
The following theorem is classical.

\begin{theorem}
Let $g$ be an element of $GL_n(\Fiel)$. The $GL_n(\Fiel)$-conjugacy class of $g$ is the $GL_n(\Fiel)$-conjugacy class  of matrix $\bigoplus_{p} J_{\nu_p, p}$ where the sum is over all irreducible factors $p$ of the minimal (or characteristic) polynomial. Here $\nu_p=\{r_1, \dots, r_k\}$ is a partition, and $\sum_p|\nu_p|\deg(p)=n$; the matrix $J_{\nu, p}=J_{r_1,p}\oplus\cdots \oplus J_{r_k,p}$ where each $J_{r_i,p}$ is a Jordan block matrix involving a companion matrix corresponding to $p$.
\end{theorem}

This theorem allows us to construct the multi-set which corresponds to the conjugacy class of $g$. It is simply $\cup_{p(t)}\{p(t)^{r_1},\dots, p(t)^{r_k}\}$ where the union is taken over all irreducible factors of the characteristic polynomial of $g$.

To classify real conjugacy classes of $GL_n(\Fiel)$ first we need to look at a Jordan block. Recall that, for a monic irreducible polynomial $p(t)$, $\tilde p(t)$ is the monic irreducible polynomial whose roots are the inverse of the roots of $p(t)$.
\begin{lemma}\label{l: irr}
Let $p(t)$ be an irreducible polynomial of degree $d$. Then $J_{r,p}^{-1}$ is conjugate in $GL_{rd}(\Fiel)$ to $J_{r,\tilde p}$, while $\zeta J_{r,p}^{-1}$ is conjugate in $GL_{rd}(\Fiel)$ to $J_{r,\breve p}$.
\end{lemma}
\begin{proof}
The theory of Jordan canonical forms tells us that $J_{r,p}^{-1}$ and $\zeta J_{r,p}^{-1}$ must be conjugate in $GL_{dr}(\Fiel)$ to Jordan block matrices. Now the characteristic polynomial of $J_{r,p}$ is $p(t)^r$ and so, by considering roots, the characteristic polynomial of $J_{r,p}^{-1}$ must be $\tilde p(t)^r$. Thus $J_{r,p}^{-1}$ is conjugate to $J_{r,\tilde p}$ as required. Similarly, the characteristic polynomial of $\zeta J_{r,p}^{-1}$ must be $\breve p(t)^r$. Thus $\zeta J_{r,p}^{-1}$ is conjugate to $J_{r,\breve p}$ as required.
\end{proof}

\begin{proposition}\label{p: glnk}
A matrix $g\in GL_n(\Fiel)$ is real if and only if $g$ is conjugate in $GL_n(\Fiel)$ to $$\left[\bigoplus_{p\neq \tilde p}(J_{\nu,p}\oplus J_{\nu,\tilde p})\right]\bigoplus \left[\bigoplus_{p=\tilde p}J_{\mu, p}\right].$$
A matrix $g$ in $GL_n(\Fiel)$ is $\zeta$-real if and only if $g$ is conjugate in $GL_n(\Fiel)$ to
$$\left[\bigoplus_{p\neq \breve p}(J_{\nu,p}\oplus J_{\nu,\breve p})\right]\bigoplus \left[\bigoplus_{p=\breve p}J_{\mu, p}\right].$$
Here $\mu$ and $\nu$ are partitions which vary with $p$.
\end{proposition}
\begin{proof}
Let $g$ be conjugate to $J_g=\oplus_p J_{\nu,p}$ and $g^{-1}$ be conjugate to $J_{g^{-1}}=\oplus_q J_{\nu,q}$. Lemma \ref{l: irr} implies that there is a matching between the Jordan blocks of $J_g$ and $J_{g^{-1}}$ which takes $J_{r,p}$ to $J_{r,q}=J_{r,\tilde p}$. Similarly there is a matching between the Jordan blocks of $J_g$ and $\zeta J_{g^{-1}}$ which takes $J_{r,p}$ to $J_{r,q}=J_{r,\breve p}$.This yields the given formulae.
\end{proof}

Proposition \ref{p: glnk} asserts that a matrix $g\in GL_n(\Fiel)$ is real (resp. $\zeta$-real) if, for any invariant factor $i$ of $g$ and for $p$ irreducible in $k[t]$, $p$ and $\tilde p$ (resp. $\breve p$) occur as factors of $i$ with the same multiplicity.

\subsection{Real conjugacy classes in $GL_n(q)$}\label{s: glnq}

We follow the notation of Macdonald in \cite{macdonald} where he gives another way to classify conjugacy classes; we will use Macdonald's method to give a criterion for an element of $GL_n(q)$ to be real. We begin by stating Macdonald's result regarding $GL_n(q)$.

\begin{theorem}\cite[1.8,1.9]{macdonald}\label{t: macdonaldconjugacy}
Let $C$ be a conjugacy class in $GL_n(q)$. Then we can associate $C$ with a sequence of polynomials $u=(u_1,u_2,\ldots)$ satisfying the following properties:
\begin{enumerate}
\item $u_i(t)=a_{n_i}t^{n_i}+\cdots+a_1t+1\in \mathbb F_q[t]$ for all $i$ with $a_{n_i}\neq 0$;
\item $\sum_{i} in_i = n$.
\end{enumerate}
This gives a one-one correspondence between conjugacy classes in $GL_n(q)$ and sequences of polynomials with the given properties.
\end{theorem}

Note that the sequence $(n_1, n_2,\dots)$ is equivalent to a partition, $\nu=1^{n_1}2^{n_2}\cdots$, of $n$; the conjugacy class $C$ described in the theorem is said to be {\it associated with} the partition $\nu$ and an element $g$ in $C$ is said to be {\it of type} $\nu$.

We need to describe how the correspondence given in Theorem \ref{t: macdonaldconjugacy} works. Let us start with a conjugacy class $C$ in $GL_n(q)$. As we described above, this can be associated with a multi-set of polynomials,
$$\{f_1(t)^{a_1}, \dots, f_r(t)^{a_r}\}$$
where, for all $i=1,\dots, r$, $f_i(t)$ is a monic irreducible polynomial in $\Fq[t]$ which is not equal to $t$, $a_i$ is a positive integer, and $\sum_{i=1}^r\textrm{deg}(f_i)a_i=n$.

Now define
$$u_i(t)=k\prod_{\{f_j(t):a_j=i\}} f_j(t).$$
Here $k\in\Fq$ is chosen so that $u_i(t)$ has constant term $1$. So $u_i(t)$ is simply the product of all irreducible polynomials in the multi-set which have associated exponent equal to $i$. That this construction gives a one-to-one correspondence with the conjugacy classes of $GL_n(q)$ is the content of Theorem \ref{t: macdonaldconjugacy}.

\begin{proposition}\label{p: correspond}
An element $g\in GL_n(q)$ is real (resp. $\zeta$-real) if and only if each of the polynomials $u_i$ in the sequence $u=(u_1,u_2,\ldots)$ (associated uniquely to the conjugacy class of $g$) are self-reciprocal (resp. $\zeta$-self-reciprocal).
\end{proposition}
\begin{proof}
Suppose that all of the $u_i(t)$ are self-reciprocal. This means that if $p(t)$ is a monic irreducible polynomial dividing $u_i(t)$, then either $p(t)=\tilde p(t)$ or $\tilde p(t)$ also divides $u_i(t)$ (with the same multiplicity). Each monic irreducible divisor of $u_i(t)$ corresponds to a Jordan block within the Jordan canonical form for $g$. Referring to Proposition \ref{p: glnk} this means that $g$ is self-reciprocal.

The converse works the same way: if $g$ is self-reciprocal then we can apply Proposition \ref{p: glnk}; thus if $p(t)$ divides $u_i(t)$ then either $p(t)=\tilde p(t)$ or else $\tilde p(t)$ also divides $u_i(t)$ (with the same multiplicity). This means that the $u_i(t)$ are self-reciprocal.

The same argument applies for the $\zeta$-real case except that we replace $\tilde p(t)$ in our argument with $\breve p(t)$.
\end{proof}

Hence to count the number of real conjugacy classes in $GL_n(q)$ we need to count sequences of polynomials $u=(u_1(t),\ldots,u_i(t),\ldots)$ such that $u_i(t)=a_it^{n_i}+\cdots+a_1t+1$ are self-reciprocal polynomials over $\Fq$ with constant term $1$ satisfying $\sum_i in_i=n$. Thus for a given partition $\nu$  we write $gl_{\nu}=\prod_{n_i>0}n_{q,n_i}$ for the number of real $GL_n(q)$-conjugacy classes of type $\nu$ in $GL_n(q)$. Then we have

\begin{theorem}\label{t: realgln}
The total number of real conjugacy classes in $GL_n(q)$ is
$$
\sum_{\{\nu:|\nu|=n\}} gl_\nu = \sum_{\{\nu:|\nu|=n\}}\prod_{n_i>0}n_{q,n_i}.
$$
Furthermore all real classes in $GL_n(q)$ are, in fact, strongly real.
\end{theorem}

We have not proved the statement about strong reality; however this follows from Corollary \ref {c: not2mod4} which we prove later on. In fact Wonenburger has proved that reality is equivalent to strong reality in $GL_n(\Fiel)$ for all fields $\Fiel$ of characteristic not $2$ \cite{wonen}.

 Note that Gow (who effectively proved Theorem \ref{t: realgln} in \cite{gow}) makes the interesting observation that the number of equivalence classes of nondegenerate bilinear forms of rank $n$ over $\Fq$ is equal to the number of real classes in $GL_n(q)$ \cite[Lemma 2.2]{gow}; hence Theorem \ref{t: realgln} also provides a formula for the number of these equivalence classes.

Note that the function $\frac{(1+t)^{(2,q-1)}}{1-qt^2}$ is a generating function for $n_{q,d}$ and so the function
$$\prod_{r=1}^{\infty} \frac{(1+t^r)^{(2,q-1)}}{1-qt^{2r}}$$
 is a generating function for the number of real conjugacy classes in $GL_n(q)$.

Finally for $q$ odd we can also write a formula for the total number of $\zeta$-real conjugacy classes in $GL_n(q)$:
$$
\sum_{\{\nu:|\nu|=n\}}\prod_{n_i>0}n_{q,n_i}\sigma_{n_i}.
$$

\section{$SL_n(q)$,  $n\not\equiv 2 \imod 4$ or $ q\not\equiv 3\imod 4$}\label{s: sl}

We count the real conjugacy classes of $SL_n(q)$ using the correspondence given by Macdonald. Note that Macdonald proves that a $GL_n(q)$ conjugacy class of type $\nu=\{\nu_1,\dots, \nu_r\}$ contained in $SL_n(q)$ is the union of $h_{\nu}$ many $SL_n(q)$ conjugacy classes where $h_\nu=(q-1, \nu_1, \dots, \nu_r)$ \cite[(3.1)]{macdonald}.

\begin{proposition}\label{p: slnu}
Let $\nu$ be a partition of $n$. Then the total number of $GL_n(q)$-real $GL_n(q)$-conjugacy classes of type $\nu$ contained in $SL_n(q)$ is
$$sl_{\nu}=\left\{ \begin{array}{ll}
\prod\limits_{n_i>0} n_{q,n_i}, & \textrm{if } q \textrm{ is even or } n_i \textrm{ is zero for } i \textrm{ odd}; \\
\frac12 \prod\limits_{n_i>0} n_{q,n_i}, & \textrm{if } q \textrm{ is odd and there exists } i \textrm{ with } in_i \textrm{ odd}; \\
f_\nu(q)\prod\limits_{i \textrm { odd}, n_i>0} q^{\frac{n_i}2-1}\prod\limits_{i \textrm { even},n_i>0}n_{q,n_i}, & \textrm{otherwise}.\end{array}\right.
$$
Here $f_\nu(q)=\frac{(q+1)^r+(q-1)^r}{2}$ where $r$ is the number of odd values of $i$ for which $n_i>0$.
\end{proposition}
\begin{proof}

As outlined in Section \ref{s: glnq} we write $\nu=1^{n_1}2^{n_2}\cdots$ with $|\nu|=n$.
Let $c$ be a conjugacy class of $GL_n(q)$ of type $\nu$ given by $u=(u_1(t),u_2(t),\ldots)$ where $u_i(t)=a_it^{n_i}+\cdots+1$.
Suppose $c$ is real and so $u_i(t)$ is self reciprocal for all $i$; this means, in particular, that $a_i=\pm 1$ for all $i$. Now $\det(1-tg)$ is that scalar multiple of the characteristic polynomial of $g$ for which the constant term equals $1$. Thus $\det(1-tg)=\prod_{i\geq 1}u_i(t)^i$ and so $\det g=(-1)^n\prod_{n_i>0} a_i^i$ \cite[(1.7)]{macdonald}. Thus a real element of $GL_n(q)$ must have determinant $\pm1$. If $p=2$ this means that all real $GL_n(q)$-conjugacy classes lie in $SL_n(q)$.

If $q$ is odd and $n_i$ is zero for $i$ odd then $n$ is even, and $\det g = \prod_{n_i>0} a_i^i = 1$; hence, again, all real $GL_n(q)$-conjugacy classes lie in $SL_n(q)$.

Now suppose that $q$ is odd and that there exists $i$ such that $in_i$ is odd. Lemma \ref{l: nqd} implies that there are $2q^{\frac{n_i-1}2}$ possibilities for $u_i(t)$; half of these will have $a_i=1$, the other half will have $a_i=-1$. Then it is easy to see that exactly half of the sequences associated to $\nu$ correspond to elements with determinant $1$.

Finally suppose that $q$ is odd and that $n_i$ is even whenever $i$ is odd, with at least one such $n_i>0$. Clearly $n$ is even so we require that $\prod_{n_i>0} a_i^i=1$. When $i$ is even, $a_i^i=1$ so we must ensure that $\prod_{i \textrm{ odd}, n_i>0} a_i^i=1$. So let us suppose, for the moment, that we are dealing with a sequence of length $r$ of self-reciprocal polynomials $(u_{i_1}(t),\dots, u_{i_r}(t))$ where $i_j$ is odd for all $j$ and $\deg u_{i_j}(t)$ is positive and even for all $j$.

Lemma \ref{l: nqd} implies that the total number of such sequences is $(q+1)^r\prod_{i \textrm { odd}, n_i>0} q^{\frac{n_i}2-1}$. Those sequences which correspond to a matrix with determinant $1$ will have an even number of leading coefficients $-1$. Counting such sequences is equivalent to counting terms in the expansion of $(q+1)^r$ in which $1$ turns up an even number of times; thus the number of such sequences is equal to $f_\nu(q)\prod_{i \textrm { odd}, n_i>0} q^{\frac{n_i}2-1}$; here $f_\nu(q)=a_r q^r+a_{r-2}q^{r-2}+a_{r-4}q^{r-4}\dots$ where $(q+1)^r = a_rq^r+a_{r-1}q^{r-1}+\dots+a_1q+a_0$. It is an easy matter to see that $f_\nu(q) = \frac{(q+1)^r+(q-1)^r}{2}$.

Now, if we return to the case where $n_i$ may be non-zero for even $i$, then it is clear that there are $\prod_{i \textrm { even},n_i>0}n_{q,n_i}$ possible polynomials corresponding to even $i$; these will make no difference to the determinant of our element, hence we obtain the given formula.

\end{proof}

Now we need to know how a real (resp. $\zeta$-real) conjugacy class of $GL_n(q)$ splits in $SL_n(q)$.

\begin{lemma}\label{l: splitup}
Let $g\in SL_n(q)$ and suppose that the conjugacy class of $g$ in $GL_n(q)$ is $C=C_1\cup\cdots\cup C_r$ where the $C_i$ are $SL_n(q)$-conjugacy classes. Take $x_i\in GL_n(q)$ such that $g^{x_i}\in C_i$ for $i=1,\dots, r$. Then $g$ is real (resp. $\zeta$-real) in $SL_n(q)$ if and only if ${g}^{x_i}$ is real (resp. $\zeta$-real) in $SL_n(q)$ for all $i$.
\end{lemma}
\begin{proof}
Suppose $tgt^{-1}=g^{-1}$ for $t\in SL_n(q)$. Then ${t}^x{g}^x(t^x)^{-1}=(g^x)^{-1}$ for $x\in GL_n(q)$. Suppose $tgt^{-1}=\zeta g^{-1}$ for $t\in SL_n(q)$. Then ${t}^x{g}^x(t^x)^{-1}=\zeta (g^x)^{-1}$ for $x\in GL_n(q)$. (Recall that, for $g\in GL_n(q)$, we write $\zeta g$ when we mean $(\zeta I)g$.)

In both cases $t^x$ is in $SL_n(q)$ since $SL_n(q)$ is normal in $GL_n(q)$. Hence we obtain our result.

\end{proof}

\begin{proposition}\label{p: not2mod4}
Suppose that $n\not\equiv 2\imod 4$ or that $q\not\equiv 3\imod 4$. Then a real (resp. $\zeta$-real) conjugacy class of $GL_n(q)$ which is contained in $SL_n(q)$ is again a union of real (resp. $\zeta$-real) conjugacy classes in $SL_n(q)$.
\end{proposition}
\begin{proof}
We generalize the proof of Wonenburger \cite{wonen}.

First consider the even characteristic situation. Let $g$ be a real element in $G=GL_n(q)$ (remember that, for even characteristic, $\zeta$-real elements do not exist) and consider the group $R_G(g)$, as defined in Section \ref{s: gl}. If $g$ is not real in $SL_n(q)$ then $R_G(g)\cap SL_n(q) \leq C_G(g)$. This implies that $R_G(g)SL_n(q)$ is a group which contains $SL_n(q)$ with even index. But this is impossible since $R_G(g)SL_n(q)$ is contained in $GL_n(q)$ which contains $SL_n(q)$ as a subgroup with odd index. We conclude that $g$ is indeed real in $SL_n(q)$.

From here on we assume that the characteristic is odd. Let $\delta_1(t),\ldots,\delta_n(t)$ be the invariant factors of $g$ in $\Fq[t]$. Suppose that $g$ is real (resp. $\zeta$-real), so each $\delta_i(t)$ is self-reciprocal (resp. $\zeta$-self-reciprocal). The $\Fq[t]$-module $V$ which we discussed in Section \ref{s: gl} decomposes as $V=\oplus_{i=1}^n V_i$, where each $V_i$ is a cyclic submodule of $V$. Write $g_i=g|_{V_i}$; then $g=\oplus_i g_i$ and the characteristic polynomial of $g_i$ is the polynomial $\delta_i(t)$.

{\bf Step 1}: We shall construct involutions in $GL(V_i)$, conjugating $g_i$ to $g_i^{-1}$ (resp. $\zeta g_i^{-1}$). If $\dim V_i=2m$ then this involution will have determinant $(-1)^m$, while if $\dim V_i=2m+1$ then we will construct two such involutions, one of determinant $1$ and the other of determinant $-1$. This construction is enough to prove the proposition except when $q\equiv 1 \pmod 4$ and $n\equiv 2 \pmod 4$; we deal with this exception in Step 2.

When $g$ is real we can write the characteristic polynomial of $g_i$ as $\chi_{g_i}(t)=(t-1)^r(t+1)^sf(t)$ where $f(\pm 1)\neq 0$ and $V_i=W_{-1}\oplus W_{1}\oplus W_0$, where $W_{-1}, W_1$ and $W_0$ are the kernels of $(g_i-1)^r, (g_i+1)^s$ and $f(g_i)$ respectively. To produce the involution $h_i$ on $V_i$ as above, it suffices to do so on each of $W_{-1}, W_1$ and $W_0$.

When $g$ is $\zeta$-real we can write $\chi_{g_i}(t)=(t^2-\zeta)^r f(t)$ where the degree of $f$ is $2m$ and $f(0)=\zeta^m$; then $V_i=W_{\zeta}\oplus W_0$, where $W_{\zeta}$ and $W_0$ are the kernels of $t^2-\zeta$ and $f(t)$ respectively. To produce the involution $h_i$ on $V_i$ as above, it suffices to do so on each of $W_{\zeta}$ and $W_0$.

It is sufficient to find a reversing involution in the following situations. Let $k$ be a cyclic linear transformation on a vector space $W$ with characteristic polynomial $\chi_{k}(t)$, one of the following three types:
\begin{enumerate}
\item $\chi_k(t)$ is self-reciprocal (resp. $\zeta$-self-reciprocal) and the degree of $\chi_{k}(t)$ is even, say $2m$. In this case we deal with the $\zeta$-real, and real, situations simultaneously. We will write the proof for the $\zeta$-real situation; the proof will hold for the real situation by replacing $\zeta$ with $1$. We assume that $\chi_k(0)=\zeta^m$.
\item $\chi_k(t)=(t-1)^{2m+1}$ or $(t+1)^{2m+1}$.
\item $\chi_k(t)=(t^2-\zeta)^m$ with $m$ odd.
\end{enumerate}
We claim that in the first case $k$ is reversed by an involution whose determinant is $(-1)^m$; in the second case there are reversing involutions with determinant $-1$, and with determinant $1$; in the third case there is an involution with determinant $-1$.

\noindent{\bf Case 1.} Since $\chi_k(t)$ is $\zeta$-self-reciprocal our assumptions imply that
$$\chi_{k}(t) = t^{2m}+a_1\zeta t^{2m-1}+a_2\zeta^2 t^{2m-2}+\dots + a_m\zeta^mt^m\dots+a_2\zeta^mt^2+a_1\zeta^mt+\zeta^m.$$
Since $W$ is cyclic, there is a vector $u\in W$ such that $\mathcal E=\{u,ku,\ldots,k^{2m-1}u\}$ is a basis of $W$.
By substituting $k^mu=y$ we get $\mathcal E = \{k^{-m}y,\ldots,y,\ldots,k^{m-1}y\}$.
Let
$$
\mathcal B=\{y,(k+\zeta k^{-1})y,\ldots,(k^{m-1}+\zeta^{m-1}k^{-m+1})y,(k-\zeta k^{-1})y,\ldots,(k^m-\zeta^mk^{-m})y\}.
$$
Then we claim that $\mathcal B$ is a basis of $W$. To see this observe first of all that, for $i=1,\dots, m-1$, $(k^{i}+\zeta^i k^{-i})y$ and $(k^i-\zeta^i k^{-i})y$ span the same $2$-dimensional subspace as $k^i y$ and $k^{-i} y$. We know also that $y$ is independent of $k^i y$ and $k^{-i} y$ hence we need only demonstrate that $(k^m-\zeta^m k^{-m})y$ is linearly independent from the rest. Were this not the case, however, we would have
\begin{eqnarray*}
&&k^m(y)=\zeta^mk^{-m}(y)+f(k^{-m+1}, k^{-m+2}, \dots, k^{m-1})(y) \\
&\implies & (k^{2m} - f(k, k^2, \dots, k^{2m-1}) - \zeta^m)y = 0
\end{eqnarray*}
where $f$ is some linear function. But this implies that
$$\chi_{k}(t) = t^{2m}-f(t^{2m-1}, t^{2m-2},\dots, t)-\zeta^m$$
which contradicts the form of $\chi_k(t)$ given above.

Now we denote the subspace generated by the first $m$ vectors of $\mathcal B$ by $P$ and that by the latter $m$ vectors by $Q$. Now observe some facts:
\begin{enumerate}
 \item $(k+\zeta k^{-1}) (k^i+\zeta^i k^{-i}) = (k^{i+1}+\zeta^{i+1} k^{-(i+1)}) + \zeta(k^{i-1}+\zeta^{i-1} k^{-(i-1)})$
\item $(k-\zeta k^{-1}) (k^i+\zeta^i k^{-i}) = (k^{i+1}-\zeta^{i+1} k^{-(i+1)}) - \zeta(k^{i-1}-\zeta^{i-1} k^{-(i-1)})$
\item $(k-\zeta k^{-1}) (k^i-\zeta^i k^{-i}) = (k^{i+1}+\zeta^{i+1} k^{-(i+1)}) - \zeta(k^{i-1}+\zeta^{i-1} k^{-(i-1)})$
\item $k^m+\zeta^mk^{-m} =-a_1(k^{m-1}+\zeta^{m-1}k^{-m+1})-a_2 (k^{m-2}+\zeta^{m-2}k^{-m+2})-\dots$
\end{enumerate}
Now (1),(2) and (4) imply that $(k+\zeta k^{-1})(P)\subseteq P$ and that $(k-\zeta k^{-1})(P)\subseteq Q$.

If we apply $k-\zeta k^{-1}$ to both sides of (4) then we find that $(k^{m+1}-\zeta^{m+1} k^{-m-1})(y)\in Q$. This, along with (2), implies that $(k+\zeta k^{-1})(Q)\subseteq Q$.

Similarly applying $k+\zeta k^{-1}$ to both sides of (4) implies that $(k^{m+1}+\zeta^{m+1} k^{-m-1})(y)\in P$. This, along with (3), implies that $(k-\zeta k^{-1})(Q)\subseteq P$.

Now let $h=1|_{P}\oplus -1|_{Q}$ and note that
$$h(k+\zeta k^{-1})h=k+\zeta k^{-1};$$
$$h(k-\zeta k^{-1})h=-(k-\zeta k^{-1}).$$
This implies that $h$ is an involution which conjugates $k$ to $\zeta k^{-1}$ and has determinant $(-1)^m$.

\noindent{\bf Case 2.} In this case, we have the characteristic polynomial $\chi_{k}(t)=(t-\epsilon)^{2m+1}$ where $\epsilon=\pm 1$.
Since $W$ is cyclic, there is a vector $u\in W$ such that $\mathcal E=\{u,ku,\ldots,k^{2m}u\}$ is a basis.
By substituting $k^mu=y$ we get $$\mathcal E = \{k^{-m}y,\ldots,y,\ldots,k^{m}y\}.$$
We consider the basis
$$
\mathcal B=\{y,(k+k^{-1})y,\ldots,(k^{m}+k^{-m})y,(k-k^{-1})y,\ldots,(k^m-k^{-m})y\}.
$$
We denote the subspace generated by the first $m+1$ vectors of $\mathcal B$ by $P$ and the latter $m$ vectors by $Q$. Examining the equation $(k-\epsilon I)^{2m+1}=0$ and applying $k\pm \epsilon I$ to both sides yields the following facts: $k+k^{-1}$ leaves $P$ as well as $Q$ invariant; also $(k-k^{-1})(P)\subseteq Q$ and $(k-k^{-1})(Q)\subset P$.
We consider $h=1|_{P}\oplus -1|_{Q}$ and $h'=-1|_P\oplus 1|_Q$.
Then $h$ and $h'$ are both involutions which conjugate $k$ to $k^{-1}$ and they have determinants $(-1)^m$ and $(-1)^{m+1}$ respectively.

\noindent{\bf Case 3.} In this case, we have the characteristic polynomial $\chi_{k}(t)=(t^2-\zeta)^m$ where $m$ is odd. Since $W$ is cyclic, there is a vector $u\in W$ such that $\mathcal E=\{u,ku,\ldots,k^{2m-1}u\}$ is a basis. By substituting $k^mu=y$ we get $$\mathcal E = \{k^{-m}y,\ldots,y,\ldots,k^{m-1}y\}.$$
We consider the basis
\begin{eqnarray*}
&\mathcal{B}=& \{y,(k-\zeta k^{-1})y,(k^2+\zeta^2 k^{-2})y, (k^3-\zeta^3 k^{-3})y, \ldots,(k^{m-1}+\zeta^{m-1}k^{-m+1})y, \\
&& (k+\zeta k^{-1})y, (k^2-\zeta^2 k^{-2})y,(k^3+\zeta^3 k^{-3})y, \ldots,(k^m+\zeta^mk^{-m})y\}.
\end{eqnarray*}
We denote the subspace generated by the first $m$ vectors of $\mathcal B$ by $P$ and the latter $m$ vectors by $Q$. This time $k-\zeta k^{-1}$ leaves $P$ as well as $Q$ invariant.
Also $(k+\zeta k^{-1})(P)\subseteq Q$ and $(k+\zeta k^{-1})(Q)\subseteq P$.
Define $h=1|_{P}\oplus -1|_{Q}$ and observe that $h$ is an involution which conjugates $k$ to $\zeta k^{-1}$ and which has determinant $(-1)^m$.

{\bf Step 2}: Write $g=\oplus_{i} g_i$, as above; the construction that we have outlined in Step~1 yields a reversing involution, $h=\oplus_i h_i$, where $h_i$ is a reversing involution for $g_i$. Whenever $n\not\equiv 2\pmod 4$ or $q$ is even, it is easy to see that the construction in Step~1 yields an involution $h$ which has determinant $1$. We must deal with the remaining situation: when $n\equiv 2 \pmod 4$ and $q\equiv 1 \pmod 4$.

Suppose that one of the submodules $V_i$ has odd dimension. Then we can choose $h_i$ to have determinant $1$ or $-1$; this ensures that we can choose $h$ to have determinant $1$ and we are done.

Suppose that all of the submodules $V_i$ have even dimension, and that the reversing involution which we constructed in Step 1, $h=\oplus_i h_i$, has determinant $-1$. Clearly one of the submodules, $V_j$ say, must have dimension $n_j\equiv 2\pmod 4$. We know that there is an $\alpha\in\Fq$ such that $\alpha^{n_j}=-1$ (see Lemma \ref{l: solxn}). Then we can define
$$h'=\left(\oplus_{i=1,i\neq j }^nh_i\right)\oplus h_j (\alpha I).$$
Clearly $h'$ is a reversing element for $g$ and $\det h'=1$ as required.
\end{proof}

\begin{corollary}\label{c: not2mod4}
If $g$ is real (resp. $\zeta$-real) in $GL_n(q)$ then $g$ is strongly real (resp. strongly $\zeta$-real) in $\langle SL_n(q), k\rangle$ for $k$ an element of determinant $-1$. If, furthermore, $n\not\equiv 2\imod 4$ or $q$ is even, then an involution $h$ exists in $SL_n(q)$ such that $hgh=g^{-1}$.
\end{corollary}
\begin{proof}
Take $g$ real (resp. $\zeta$-real)in $GL_n(q)$. Note that in the proof of Proposition \ref{p: not2mod4}, we did not use the fact that $g$ is contained in $SL_n(q)$ and, in all cases, we found a reversing element of determinant $1$ or $-1$. In fact, in the odd characteristic case, we always found a reversing element which was an involution; this proves the first statement when the characteristic is odd.

For $n\not\equiv 2 \pmod 4$ and $q$ odd, we were able to do better: we found an involution in $SL_n(q)$ which reverses $g$; this proves the second statement when the characteristic is odd. The only thing left to prove is the following: if $g$ is real in $GL_n(q)$ with $q$ even, then $g$ is strongly real in $GL_n(q)$. (This statement is strong enough because, when $q$ is even, there are no $\zeta$-real elements, and all involutions of $GL_n(q)$ lie in $SL_n(q)$.)

Take $g$ real in $GL_n(q)$ with $q$ even. We refer to Proposition \ref{p: glnk} and write $g$ as a block matrix as follows:
$$\left[\bigoplus_{p\neq \tilde p}(J_{\lambda,p}\oplus J_{\lambda,\tilde p})\right]\bigoplus \left[\bigoplus_{p=\tilde p}J_{\lambda,p}\right].$$
We will consider two cases:
\begin{enumerate}
 \item $g= J_{r,p}\oplus J_{r,\tilde p}$
\item $g= J_{r,p}$ and $p=\tilde p$.
\end{enumerate}
In both cases $p$ is an irreducible polynomial. Clearly if we can find an involution which reverses $g$ in both of these cases, then we can build an involution which reverses $g$ in general.

Consider the first case. Then $g$ is conjugate to a matrix of form
$$g_1=\left(\begin{array}{cc}
B & 0 \\
0 & B^{-1}
        \end{array}\right),$$
where $B$ is a square $d$ by $d$ matrix for $d$ the degree of $p$. Now observe that $hg_1h^{-1}=g_1^{-1}$ where
$$h=\left(\begin{array}{cc}
0 & I \\
I & 0
        \end{array}\right).$$
Since $h^2=1$, we conclude that $g_1$ is strongly reversible and thus so is $g$.

Now consider the second case. If $p=t+1$ then $g$ is an involution and we are done. Thus we assume that $p$ has even degree, $d$. We write $g=g_s g_u$, the Jordan decomposition into semi-simple and unipotent elements. Since $g$ is real, Lemma \ref{l: centralizer} implies that $g_s$ is real. Now $C_{GL_n(q)}(g_s)\cong GL_r(q^d)$; the element which reverses $g$ we denote by $x$; it acts as an involutory field automorphism on $GL_r(q^d)$.

We can think of $g$ as lying inside $GL_r(q^d)$; then $g$ is conjugate to an element $g_1$ where
$$g_1=\left(\begin{array}{cccc}
\alpha & 1 && \\
&\ddots & \ddots & \\
&&\ddots&1 \\
 &&& \alpha \\
        \end{array}\right),
\
g_1^{-1}=\left(\begin{array}{ccccc}
\alpha^{-1} & -\alpha^{-2} & \alpha^{-3} && \\
&\ddots & \ddots & \ddots & \\
&&\ddots& \ddots & \alpha^{-3} \\
&&&\ddots & -\alpha^{-2} \\
 &&&& \alpha^{-1} \\
        \end{array}\right)
$$
where $\alpha$ is an element of $\mathbb{F}_{q^d}$ which satisfies $\alpha^{q^{\frac{d}2}} = \alpha^{-1}$. Then $g_1$ is reversed by an element $h \sigma$ where
$$h=\left(\begin{array}{cccccc}
1 & 0& 0 & \cdots && \\
& -\alpha^2 & \alpha^3 & -\alpha^4 & \cdots& \\
&& \alpha^4 & -2\alpha^5 & 3\alpha^6 & \cdots \\
&&& -\alpha^6 & 3\alpha^7 & \cdots \\
&&&& \ddots &
        \end{array}\right), \ \textrm{i.e. }
h_{ij}=\left\{ \begin{array}{ll}
(-1)^{j+1}\left(\genfrac{}{}{0pt}{}{j-2}{i-2}\right)\alpha^{i+j-2}, & j\geq i; \\
0, & j<i. \end{array}\right.
$$
Here we define $\left(\genfrac{}{}{0pt}{}{k}{0}\right)=1$ and $\left(\genfrac{}{}{0pt}{}{k}{-1}\right)=-1$ for $k$ a positive integer (of course, since the characteristic is even, $h_{ij}=0$ for many $i,j$). Now $(hx)^2=h h^x$ acts trivially on $GL_r(q^d)$. But this means that $(hx)^2\in Z(GL_r(q^d))$. Since $Z(GL_r(q^d))$ has odd order we conclude that $hx$ can be chosen so that $(hx)^2=1$. Hence $g_1$ is strongly real and thus so is $g$.

\end{proof}

We summarise our results for real elements with the following theorem:

\begin{theorem}\label{t: countinGL}
Suppose that $n\not\equiv 2\imod 4$ or $q\not\equiv 3\imod 4$. Then the total number of real conjugacy classes in $SL_n(q)$ is equal to:
$$\sum_{|\nu|=n}h_\nu sl_{\nu} $$
where $\nu=\{\nu_1, \dots, \nu_r\}$, $h_\nu=(q-1, \nu_1,\dots, \nu_r)$ and the value of $sl_{\nu}$ is given in Proposition \ref{p: slnu}. Furthermore, provided $n\not\equiv 2\imod 4$ or $q$ is even, this is the same as the total number of strongly real conjugacy classes in $SL_n(q)$.
\end{theorem}

\section{$SL_n(q)$,  $n\equiv 2 \imod 4$  and  $q\equiv 3\imod 4$}\label{s: sl2}

In this section we assume that $n\equiv 2\imod 4$ and $q\equiv 3\imod 4$, and we take $\zeta=-1$, a non-square. The formula given in Theorem \ref{t: countinGL} gives the number of conjugacy classes in $SL_n(q)$ which are real in $GL_n(q)$. Thus, in order to count the number of real conjugacy classes in $SL_n(q),$ we take this formula and count how many $GL_n(q)$-real conjugacy classes in $SL_n(q)$ fail to be real in $SL_n(q)$. Our analysis is based on the following:

\begin{lemma}\label{l: square}
Let $g\in SL_n(q)$. Then $g$ is real (resp. $\zeta$-real) in $SL_n(q)$ if and only if there exists a reversing element $h\in GL_n(q)$ such that $\det h$ is a square in $\Fq$.
\end{lemma}
\begin{proof}
Clearly if $g$ is real (resp. $\zeta$-real) in $SL_n(q)$ then there exists a reversing element $h$ in $SL_n(q)$ and $\det h$ is a square in $\Fq$.

Now for the converse: suppose that there exists a reversing element $h\in GL_n(q)$ such that $\det h$ is a square. We note that for a positive integer $c$ we have $\det h^c=(\det h)^c$ and, if $c$ is odd, then $h^cgh^{-c}=hgh^{-1}$ and so $h^c$ is a reversing element. Since $\det h^c=1$ for some odd integer $c$, we are done.
\end{proof}

Now let $g$ be a real (resp. $\zeta$-real) element in $GL_n(q)$ with reversing element $h\in GL_n(q)$. In view of Lemma \ref{l: square} we will be interested in determining whether $\det h$ is a square or a non-square in $\mathbb{F}_q$. We will make use of two commutative diagrams:

$$
\begin{CD}
GL_a(\mathbb{F}_q) @>{i}>> GL_a(\overline{\mathbb{F}_q})\\
@V{\det}VV                                    @VV{\det}V\\
\mathbb{F}_q^* @>>{i}> \overline{\mathbb{F}_q}^*
\end{CD} \quad \quad \quad
\begin{CD}
GL_a(q^d) @>{i}>> GL_{ad}(q)\\
@V{\det}VV                                    @VV{\det}V\\
\mathbb{F}_{q^d}^* @>>N> {\mathbb{F}_q}^*
\end{CD}
$$

The map $i$ denotes a natural inclusion map. The map $N$ is the norm map defined as follows:
$$N: C_{q^d-1}\to C_{q-1}, \quad x \mapsto x^{q^{d-1}+\dots+q+1}.$$
The commutativity of the first diagram is obvious. The commutativity of the second is explained by the following fact: $N\circ \det: GL_a(q^d)\to \Fq^*$, when viewed as a map from a subgroup of $GL_{ad}(q)$, is multilinear, alternating (on columns) and satisfies $(N\circ\det)(I)=1$; in other words it is a determinant and so must coincide with $\det\circ i: GL_a(q^d)\to\Fq^*$ by \cite[Proposition 4.6, p.514]{lang}.

Thanks to the given inclusion maps there are several determinant maps applicable to any given matrix. In what follows we write $\det_{|k|}$ to specify the field $k$ in which our image lies. Note that the form of $N$ implies that, for $y\in GL_a(q^d)$, $\det_{q}(y)$ is square if and only if $\det_{{q^d}}(y)$ is square.

We will build a reversing element for $g$, in a similar way to the proof of Corollary \ref{c: not2mod4}, by considering three basic cases. Let $V$ be the $\Fq[t]$-module associated with $V$ (see Section \ref{s: glnk}) and  let $\chi_g(t)$ be the characteristic polynomial of $g$.
 \begin{enumerate}
\item $V$ is a cyclic module for $g$ with $\chi_g(t)=p(t)=(t\pm 1)^a$.
\item $V$ is a cyclic module for $g$ with $\chi_g(t)=p(t)^a$, where $p(t)$ is a self-reciprocal (or $\zeta$-self-reciprocal) irreducible polynomial of even degree $d$.
\item $V=W_p\oplus W_q$, a module for $g$, such that $W_p$ and $W_q$ are cyclic modules with characteristic polynomials $p(t)^a$ and $\tilde p(t)^a$ (or $p(t)^a$ and $\breve p(t)^a$) such that $p(t)$ is irreducible.
\end{enumerate}

\begin{lemma}\label{l: cyclic}
Let $V$ be a cyclic module associated with $g$ such that $\chi_g(t)=(t\pm 1)^a$. Suppose that $h$ reverses $g$.
Then,
\begin{enumerate}
\item if $a$ is odd, then $\det h$ may be a square or a non-square in $GL_a(q)$.
\item if $a\equiv 0 \imod 4$, then $\det h$ is a square in $GL_a(q)$.
\item if $a\equiv 2 \imod 4$, then $\det h$ is a non-square in $GL_a(q)$.
\end{enumerate}
\end{lemma}
\begin{proof}
Write $g$ in upper triangular form. Take $h$ such that $h^{-1}gh=g^{-1}$; then $h$ must be upper triangular with diagonal $(\alpha,-\alpha,\alpha,-\alpha,\cdots)$. The result follows by considering the determinant.
\end{proof}

\begin{lemma}\label{l: cyclic2}
Let $V$ be a cyclic module associated with $g$ such that $\chi_g(t)=p(t)^a$ where $p(t)$ is an irreducible polynomial of even degree $d$. Suppose that $h$ reverses $g$. If $a$ is odd, then $\det h$ may be a square or non-square in $GL_{ad}(q)$ while, if $a$ is even, then $\det h$ is a square in $GL_{ad}(q)$.
\end{lemma}
\begin{proof}
Write $g=g_sg_u$ for the Jordan decomposition of $g$. In $GL_{d}(q)$, $g_s$ is centralized by $GL_a(q^d)$ and so $g_s$ can only be reversed by a field automorphism of $GL_a(q^d)$. Let $x$ be such a field automorphism; then $h=yx$ where $y$ is an element of $GL_a(q^d)$ which satisfies $(g_u)^{-1} = (yx) (g_u) (yx)^{-1}$ (see Lemma \ref{l: centralizer}).

If $a$ is odd then there is a central element in $GL_a(q^d)$ which has determinant a non-square in $\mathbb{F}_q^{d}$. Hence $\det h$ may be a square or a non-square in $GL_{ad}(q)$.

Assume that $a$ is even. Proposition \ref{p: not2mod4} implies that there exists a reversing element $h$ in $GL_{ad}(q)$ such that $\det h$ is a square. Any other reversing element in $GL_{ad}(q)$ must equal $hz$ where $z\in C_{GL_{ad}(q)}(g)$. If $z$ centralizes $g$ then it centralizes $g_s$ and we conclude that $z$ must be an element in $GL_a(q^d)$ which centralizes $g_u$. We write $g_u$ as an element of $GL_a(q^d)$:
$$g=\left(\begin{array}{cccc}
1 & \alpha &&  \\
&\ddots &\ddots& \\
&&\ddots & \alpha \\
&&& 1
  \end{array}\right),$$
where $\alpha\in \mathbb{F}_{q^d}$. Then $z$ must have form
$$z=\left(\begin{array}{cccc}
\beta_1 & \beta_2&&  \\
&\ddots &\ddots& \\
&&\ddots & \beta_2 \\
&&& \beta_1
  \end{array}\right).$$
Hence, since $a$ is even, $\det_{q^d} z$ is a square. Thus any reversing element for $g$ in $GL_{ad}(q)$ has determinant a square.
\end{proof}

\begin{lemma}\label{l: sltwocycles}
Let $V=W_p\oplus W_q$ be the module associated with $g$. Suppose that $W_p$ and $W_q$ are cyclic modules with characteristic polynomials $p(t)^a$ and $\tilde p(t)^a$ (or $\breve p(t)^a$). Suppose that $p(t)$ is irreducible of degree $d$, and that $h$ reverses $g$. If $a$ is odd then $\det(h)$ may be square or non-square while, if $a$ is even, then $\det(h)$ is square.
\end{lemma}
\begin{proof}
We know that $g$ is conjugate in $GL_{2ad}(q)$ either to $g_1=\left(\begin{array}{cc}
B & 0 \\
0 & B^{-1}  \end{array}\right)$ or to $g_2=\left(\begin{array}{cc}
B & 0 \\
0 & \zeta B^{-1}  \end{array}\right)$. By Lemma \ref{l: splitup} we can assume that $g$ is equal to $g_1$ or $g_2$. If $a$ is odd then $GL_a(q^d)$ contains a central element with non-square determinant in $\Fqd$; since $g$ is central in a group isomorphic to $GL_a(q^d)\times GL_a(q^d)$ we conclude that $\det(h)$ may be square or non-square. 

Suppose that $a$ is even.  If $h$ preserves blocks then $p=\tilde p$ (or $p=\breve p$) and Lemma \ref{l: cyclic2} implies that $\det(h)$ is a square. Otherwise $h$ reverses blocks and $h=\left(\begin{array}{cc}
0 & X \\
Y & 0  \end{array}\right)$ for some $X$ and $Y$ in $GL_{ad}(q)$. Thus 
$$\left(\begin{array}{cc}
0 & X \\
Y & 0  \end{array}\right)\left(\begin{array}{cc}
B & 0 \\
0 & C  \end{array}\right)\left(\begin{array}{cc}
0 & X \\
Y & 0  \end{array}\right)^{-1}= \left(\begin{array}{cc}
C & 0 \\
0 & B  \end{array}\right),$$
where $C=B$ in the real case, and $C=\zeta B^{-1}$ in the $\zeta$-real case.
This implies that $XCX^{-1}=B^{-1}$ and $YBY^{-1}=B$; now $C_{GL_{ad}(q)}(B)=C_{GL_{ad}(q)}(C)$, hence we need to examine the centralizer of $B$ in $GL_{ad}(q)$.

Clearly $\det_{q} h=(-1)^{ad}(\det X)(\det Y) = (\det X)(\det Y)$ since $a$ is even. Then $X$ must lie in the centralizer of $g_u|_{W_q}$ and we have seen the form of such a centralizer in Lemma \ref{l: cyclic2}; we know that the determinant is always a square in $\Fqd$ hence is a square in $\Fq$.
\end{proof}

This concludes our treatment of the three specific cases. We use the lemmas to build up the picture for general $V$.

\begin{proposition}\label{p: n2mod4}
Suppose that $g$ lies in $SL_n(q)$ with $n\equiv 2 \pmod 4$ and $q\equiv 3\pmod 4$. Suppose that the corresponding module, $V$, for $g$ splits into $r$ cyclic submodules with corresponding polynomials $p_1(t)^{a_1}, \dots, p_r(t)^{a_r}$. If $g$ is real in $GL_n(q)$, then $g$ is real in $SL_n(q)$ if and only if $a_i$ is odd for some $i$. If $g$ is $\zeta$-real in $GL_n(q)$ then $g$ is $\zeta$-real in $SL_n(q)$.
\end{proposition}
\begin{proof}
Suppose first that $g$ is real in $GL_n(q)$. Let $h$ be any element of $GL_n(q)$ such that $hgh^{-1}=g^{-1}$. We can break $V$ up into submodules of the three types listed above (call these $h$-{\it minimal}). If $a_i$ is odd for some $i$ then consider $W$, the $h$-minimal submodule corresponding to $a_i$. Lemmas \ref{l: cyclic} to \ref{l: sltwocycles} imply that $g|_{W}$ is conjugate to its inverse by elements of determinant $+1$ and $-1$ in $GL(W)$; this property will then hold for $g$.

Conversely if all of the $a_i$ are even then our above calculations show that, restricted to any $h$-minimal submodule, $W$, $g|_{W}$ is conjugate to its inverse by an element of determinant $+1$ or $-1$ in $GL(W)$, but not both. In fact, of the three types listed above, the only time this determinant is $-1$ is when $W$ is cyclic and $(p_i(t))^{a_i}=(t\pm1)^{a_i}$ with $a_i\equiv 2\pmod 4$ (this is also the only time the dimension of $W$ is not equivalent to $0\pmod 4$).

Since $n\equiv 2\pmod 4$ there will be an odd number of these $-1$ $h$-minimal submodules, thereby ensuring that $g$ is only conjugate to its inverse by an element of determinant $-1$.

Now suppose that $g$ is $\zeta$-real and all the $a_i$ are even. This implies that the dimension of all the $h$-minimal submodules is divisible by $4$. Since $n\equiv 2\pmod4$ this is a contradiction. Thus $a_i$ is odd for some $i$. Let $W$ be the $h$-minimal submodule corresponding to $a_i$. Once again, Lemmas \ref{l: cyclic} to \ref{l: sltwocycles} imply that $g|_{W}$ is conjugate to its inverse by elements of determinant $+1$ and $-1$ in $GL(W)$; this property will then hold for $g$.
\end{proof}

\begin{corollary}
Suppose that $g\in SL_n(q)$ is of type $\nu=1^{n_1}2^{n_2}\cdots$ and is real in $GL_n(q)$. Then $g$ is real in $SL_n(q)$ if and only if $n_i>0$ for some odd $i$.
\end{corollary}
\begin{proof}
We simply need to convert the criterion given by Proposition \ref{p: n2mod4} into the language of Macdonald, as described in Section \ref{s: glnq}.
\end{proof}

We can use this corollary to count the real classes in $SL_n(q)$ as follows:

\begin{theorem}\label{t: countinSL}
Suppose that $n\equiv 2\imod 4$ and $q\equiv 3\imod 4$. Then the total number of real conjugacy classes in $SL_n(q)$ is equal to:
$$\sum_{|\nu|=n}h_{\nu}sl_{\nu}-\sum_{|\mu|=n}h_{\mu}sl_{\mu} $$
where $\nu=(1^{n_1}2^{n_2}\ldots)$, $\mu=(2^{d_2}4^{d_4}\ldots)$ and the values of $sl_{\nu}$ and $sl_{\mu}$ are given in Proposition \ref{p: slnu}.
\end{theorem}

\section{Strongly real conjugacy classes in $SL_n(q)$}\label{s: slstrong}

\begin{theorem}\label{t: slstrongly}
Let $g$ be an element of $SL_n(q)$ which is real in $GL_n(q)$. Let $g$ be of type $\nu=1^{n_1}2^{n_2}\cdots$, with associated self-reciprocal polynomials $u_i$ of degree $n_i$.
Then,
\begin{enumerate}
\item if $n\not\equiv 2\imod 4$ or if $q$ is even then $g$ is real as well as strongly real in $SL_n(q)$.
\item if $n\equiv 2\imod 4$ and $q$ is odd then, $g$ is strongly real in $SL_n(q)$ if and only if there is an odd $i$ for which $\pm1$ appears as a root of $u_i(t)$.
\end{enumerate}
\end{theorem}
\begin{proof}
The first statement follows directly from Corollary \ref{c: not2mod4}.  Now suppose that $n\equiv 2\imod 4$ and $q$ is odd.

Take $g$ a strongly real element in $GL_n(q)$. We use the same notation as the previous section except that this time we require that $h^2=1$. Let $W$ be a $h$-minimal submodule of $V$; $W$ will have one of the same three types as before.

Suppose that $W=W_p\oplus W_q$ where $W_p$ and $W_q$ are cyclic with corresponding characteristic polynomials $p(t)^a$ and $\tilde p(t)^a$ such that $p(t)$ and $\tilde p(t)$ are distinct irreducible polynomials of degree $d$. Then $h|_W=\left(\begin{array}{cc}
0 & X \\
Y & 0  \end{array}\right)$ for some $X$ and $Y$ in $GL_{ad}(q)$ and $\det h|_W=(-1)^{ad}(\det X)(\det Y )$. But, since $h^2=1$, we have $X=Y^{-1}$ and $\det{h|_W}=(-1)^{ad}$.

Suppose that $W$ is cyclic and $d>1$. Write the corresponding characteristic polynomial as $p(t)^a$ where $p(t)$ has even degree $d$. We can consider $W\otimes_{\Fq} \overline{\mathbb{F}_q}$. Then $g|_W$ is clearly real in $GL(V\otimes_{\Fq} \overline{\Fq})$ but is reducible. In fact there are pairings of blocks. As we have already seen this implies that we have determinant $(-1)^b$ over these pairings, where $b$ is the size of the block. Thus $h$ must have determinant $(-1)^{\frac{ad}2}$.

Suppose that $W$ is cyclic and $p(t)=t\pm1$. It is easy to check that
\begin{enumerate}
\item if $a$ is odd then $\det  h|_W$ may be $1$ or $-1$.
\item if $a\equiv 0 \imod 4$ then $\det h|_W=1$.
\item if $a\equiv 2 \imod 4$ then $\det h|_W=-1$.
\end{enumerate}

We have now treated the three types. If $V$ does not have a cyclic submodule corresponding to a polynomial $(t\pm1)^a$ where $a$ is odd, then $\det h|_W = (-1)^\frac{\dim W}2$ and so $\det h = (-1)^\frac{n}2 = -1$; in particular $g$ is not strongly real in $SL_n(q)$. On the other hand if $V$ does have a cyclic submodule corresponding to a polynomial $(t\pm1)^a$ where $a$ is odd, then we can choose $h$ to have $\det h=1$.

The proof is completed once we observe that $V$ has a cyclic submodule corresponding to a polynomial $(t\pm1)^a$ where $a$ is odd if and only if there is an odd $i$ for which $\pm 1$ is a root of $u_i(t)$.
\end{proof}

\section {$PGL_n(q)$}\label{s: pgl}

First, some notation: consider two groups, $X$ and $Y$, such that $Y\leq Z(X)$. We say that an element $h\in X/Y$ {\it lifts} to an element $g$ in $X$ (or, equivalently, $g$ {\it projects onto} $h$) if $h=gY$. Now suppose that $W<X$. We say that $g$ {\it projects into} $W/Y$ if there exists $y\in Y$ such that $gy\in W$. Recall also that, for $g\in GL_n(q)$ and $\eta\in\Fq^*$, we will abuse notation and write $\eta g$ for $g(\eta I)$.

\subsection{Conjugacy in $PGL_n(q)$}\label{s: conjugacy in pgl}

Set $Z=Z(GL_n(q))$; then $PGL_n(q)=GL_n(q)/Z$; our first job is to understand how conjugacy in $PGL_n(q)$ works. Let $g$ be an element of $C$, a conjugacy class of $GL_n(q)$, represented by $u=(u_1(t),u_2(t), \dots)$ corresponding to a partition $\nu=1^{n_1}2^{n_2}\cdots$. Macdonald asserts then the conjugacy class of $\eta g$ is represented by $(u_1(\eta t), u_2(\eta t),\dots)$ \cite[p.30]{macdonald}. Thus all elements in $GL_n(q)$ which project onto an element $gZ\in PGL_n(q)$ are of type $\nu$; we therefore refer to $gZ$ as being {\it of type} $\nu$.

Suppose that $g$ projects onto $gZ$ which is real in $PGL_n(q)$. Then we want to calculate how many real and $\zeta$-real elements in $GL_n(q)$ project onto $gZ$. We define two sequences of self-reciprocal polynomials (resp. $\zeta$-self-reciprocal polynomials), $u=(u_1(t),u_2(t), \dots)$ and $v=(v_1(t),v_2(t),\dots),$ to be {\it equivalent} if, for some $\eta\in\Fq^*$, $v_i(t)=u_i(\eta t)$ for all $i$.

To understand what this means for reality in $PGL_n(q)$ we need to return to the study of self-reciprocal and $\zeta$-self-reciprocal polynomials that we started in Section \ref{s: poly}. 

\subsection{An action of $\mathbb F_q^*$ on Polynomials}\label{s: action}
We define an action of $\mathbb F_q^*$ on the set of degree $n$ polynomials by $\eta.f(t)=f(\eta t)$ for $\eta\in \mathbb F_q^*$. Recall the definition of sets $T_n$ and $S_n$ given in Section \ref{s: poly}. We are interested in classifying the orbits of $\Fq^*$ intersected with $T_n$ and $S_n$. That is to say, we wish to determine the size of the sets
$$[f]_T=\{f(\eta t)\in T_n\mid \eta\in \mathbb F_q^*\}, \textrm{ and } [f]_S=\{f(\eta t)\in S_n\mid \eta\in \mathbb F_q^*\},$$
for a degree $n$ polynomial $f$ in $\Fq[t]$.

In what follows we write $|k|_2$ for the largest power of $2$ which divides an integer $k$. That is $|k|_2=2^r$ where $2^r$ divides $k$, but $2^{r+1}$ does not. We begin with an easy arithmetical lemma:

\begin{lemma}\label{l: solxn}
Let $\Fq$ be a finite field with $q$ odd. Then there exists $\alpha\in \Fq^*$ with $\alpha^n=-1$ if and only if $|n|_2<|q-1|_2$.
\end{lemma}
\begin{proof}
Recall that $\Fq^*$ is a cyclic group of order $q-1$. Consider the homomorphism $\Fq^*\to \Fq^*, x \mapsto x^n$. The kernel of this homomorphism is the subgroup of $\Fq^*$ given by $\{x\in \Fq^* \mid x^{(n,q-1)}=1\}$. Hence the image of the homomorphism has $r=\frac{(q-1)}{(n,q-1)}$ elements. This image is the subgroup of $\Fq^*$ given by $\{x\in \Fq^*\mid x^r=1\}$. Now $-1$ is contained in this image if and only if $r$ is even or, equivalently, if and only if $|n|_2 < |q-1|_2$.
\end{proof}

\begin{lemma}\label{l: lift}
If $q$ is even then $[f]_T$ and $[f]_S$ contain at most one element. If $q$ is odd then $[f]_T$ and $[f]_S$ contain at most two elements.
\end{lemma}
\begin{proof}
Write $X$ to mean either $T$ or $S$. Take $f(t)\in X_n$ such that $f(\eta t)\in X_n$; then $\eta^{n}=\pm 1$. Since $f(\eta t)\in X_n$, for any coefficient $a_k\neq 0$ of $f(t)$, we must have $a_k\eta^{n-k}=\pm a_k\eta^k$. Thus, if the order of $\eta$ is denoted $e$, then $e|2k$. 

If $q$ is even this implies that $e|k$ and $f(t)\in\Fq[t^e]$; thus $f(\eta t)=f(t)$ as required. Suppose that $q$ is odd. If $e$ is odd then $e|k$ and, again, $f(\eta t)=f(t)$. So assume that $e$ is even. We must have $f(t)\in \Fq[t^\frac{e}{2}]$ and so $f(\eta^2 t)=f(t)$. 

Now suppose that $f(\epsilon t)\in X_n$ for some $\epsilon \in \Fq^*$. To avoid redundancy, let us assume that $f(\eta t)\neq f(t)\neq f(\epsilon t)$. Clearly $|d|_2=|e|_2$, where $d$ is the order of $\epsilon$, and $f(t)\in \Fq[t^\frac{e}{2}]\backslash \Fq[t^e]$. But then $\epsilon$ and $\eta$ are doing the same thing to coefficients of $f(t)$, namely swapping the sign of coefficients of $t^{\frac{ae}2}$ for odd $a$. Thus $f(\eta t)=f(\epsilon t)$ as required.
\end{proof}

Note that, when $q$ is odd, we have effectively shown that if $f\in X_n$ and $[f]_X$ has two elements then $[f]_X=\{f(t), f(\eta t)\}$ and the order of $\eta$ is a power of $2$. We continue our analysis for $q$ odd.

\begin{lemma}\label{l: lift2}
Let $q$ be odd.
\begin{enumerate}
\item If $n$ is odd then $S_n$ is empty and, for $f(t)\in T_n,$ we have $|[f]_T|=2$.
\item If $n$ is even then,
\begin{enumerate}
\item if $f(t)\not\in \Fq[t^{|q-1|_2}]$ and $f(t)\in T_n$ (resp. $f(t)\in S_n$) we have $|[f]_T|=2$ (resp. $|[f]_S|=2$). Moreover if $[f]_T$ is non-empty (resp. $[f]_S$ is non-empty), then $[f]_S$ is empty (resp. $[f]_T$ is empty).
\item if $f(t)\in \Fq[t^{|q-1|_2}]$ and $f(t)\in T_n$ (or $f(t)\in S_n$) then $|[f]_S|=|[f]_T|=1$. 
\end{enumerate}
\end{enumerate}
\end{lemma}
\begin{proof}
First we consider when $n$ is odd. There are no $\zeta$-self-reciprocal polynomials in this case, hence $S_n$ is empty. What is more $T_n=F_n\cup G_n$ and we have a bijection between $F_n$ and $G_n$ given by $f(t)\mapsto f(-t)$. Thus, for $f\in T_n$, $[f]_T=\{f(t), f(-t)\}$ as required.

Let us consider $n$ even now, so $|n|_2\geq 2$. Again write $X$ for either $T$ or $S$ and suppose that $f(t)\in X_n$. We know that $[f]_X=\{f(t), f(\eta t)\}$ where $\eta$ has order a power of $2$. If $f(t)\in \Fq[t^{|q-1|_2}]$ then this implies that $f(\eta t)=f(t)$ and so $|[f]_X|=1$. On the other hand if $f(t)\not\in \Fq[t^{|q-1|_2}]$ then let $e$ be the smallest power of $2$ such that $f(t)\not\in\Fq[t^e]$. Take $\eta$ of order $e$ and it is easy to check that $f(\eta t)\in X_n$ and $f(\eta t)\neq f(t)$.

Suppose that $f(t)\in S_n$ (i.e. $f(t)$ is $\zeta$-self-reciprocal) and $f(t)\not\in\Fq[t^{|q-1|_2}]$. We check if $f(\lambda t)$ is self-reciprocal for some $\lambda\in\Fq$. Then this implies that, whenever $a_k\neq 0$,
$$a_k\frac{\lambda^{n-k}}{\zeta^{\frac{n}2-k}} = \pm a_k\lambda^k.$$
This implies that $\lambda^{n-k}=\pm \zeta^{\frac{n}2}.$ Setting $k=0$ we find that $|q-1|_2$ divides $|n|_2$ and so $\lambda ^{2k}=\pm \zeta ^k$ for $k>0$. Since $\zeta$ is a non-square this is impossible. Thus, if $[f]_S$ is non-empty, then $[f]_T$ is empty. A similar argument shows that, if $[f]_T$ is non-empty, then $[f]_S$ is empty.

Finally consider the case when $f(t)\in\Fq[t^{|q-1|_2}]\cap T_n$. Let $\eta$ be a non-square of order $|q-1|_2$. Set $\kappa$ to be an element in $\Fq^*$ which satisfies $\kappa^{2} = \frac{\eta}{\zeta}$. Then one can check that
$$f(\kappa t) = \pm (\kappa t)^{n}+a_{|q-1|_2} (\kappa t)^{n-|q-1|_2}\pm a_{2|q-1|_2} (\kappa t)^{n-2|q-1|_2}+\dots \pm a_{|q-1|_2}(\kappa t)^{|q-1|_2}+1$$
lies in $S_n$. Similarly if $f(t)\in\Fq[t^{|q-1|_2}]\cap S_n$ then $f(\theta t)$ lies in $T_n$ where $\theta$ is an element of $\Fq^*$ which satisfies $\theta^2 = \eta \zeta$. Thus, in both cases, $|[f]_S|=|[f]_T|=1$.
\end{proof}

\subsection{Real classes in $PGL_n(q)$}

We wish to use our results concerning reality in $GL_n(q)$ to classify reality in $PGL_n(q)$. In general, when converting our results from a group $X$ to $X/Y$ where $Y\leq Z(X)$, we are faced with the following problem: if $g$ is real (resp. strongly real) in $X$ then $gY$ is real (resp. strongly real) in $X/Y$, however the converse does not hold.

When $X=GL_n(q)$ we are able to give a partial converse. Write $Z$ for the centre of $GL_n(q)$ and recall that $\zeta$ is a fixed non-square in $\Fq$.

\begin{lemma}\label{l: projreal}
Suppose that $gZ$ is real in $PGL_n(q)$. Then $gZ$ lifts to a real or a $\zeta$-real element in $GL_n(q)$.
\end{lemma}
\begin{proof}
Clearly $gZ$ lifts to an element $g$ which is conjugate in $GL_n(q)$ to $g^{-1} (\frac1{\eta} I)$ for some $\eta\in\Fq^*$, i.e. there exists $h\in GL_n(q)$ such that $hgh^{-1}=\eta^{-1}g^{-1}$.
If $\eta$ is a square, $\eta=\lambda^2$ say, we get  $h(\lambda g)h^{-1}=\lambda^{-1}g^{-1}$. That is, $gZ$ lifts to a real element $\lambda g$ in $GL_n(q)$.
When $\eta$ is not a  square, we write $\eta=\zeta^{-1}\lambda^2$ and we have $h(\lambda g)h^{-1}=\zeta(\lambda g)^{-1}$. In this case $gZ$ lifts to a $\zeta$-real element $\lambda g$ in $GL_n(q)$.
\end{proof}

Now the number of real classes in $PGL_n(q)$ is equal to
$$\sum_{|\nu|=n}pgl_\nu$$
where $pgl_\nu$ is the number of real conjugacy classes in $PGL_n(q)$ of type $\nu$. Lemma \ref{l: projreal} implies that $pgl_\nu$ is equal to the number of equivalence classes in the set of sequences of self-reciprocal and $\zeta$-self-reciprocal polynomials associated with $\nu$ (or, equivalently, the number of equivalence classes of real and $\zeta$-real conjugacy classes in $GL_n(q)$ of type $\nu$). 

For the rest of this section we calculate $pgl_\nu$ for different $\nu, n$ and $q$.

\subsubsection{$q$ is even}
There are no $\zeta$-real conjugacy classes in $GL_n(q)$ in this case. What is more, Lemma \ref{l: lift} implies that all equivalence classes of real conjugacy classes in $GL_n(q)$ are of size $1$; hence $pgl_\nu=gl_\nu=\prod_{n_i>0}n_{q,n_i}$.

\subsubsection{$q$ is odd}

\begin{lemma}
If $C$ is a real (resp. $\zeta$-real) class in $GL_n(q)$ then $C$ is equivalent to at most one other real (resp. $\zeta$-real) class in $GL_n(q)$
\end{lemma}
\begin{proof}
The proof is very similar to the proof of Lemma \ref{l: lift}.  Suppose that $C$ is real (resp. $\zeta$-real) and corresponds to the sequence $u=(u_1(t), u_2(t), \dots)$. Consider conjugacy classes $C_\eta$ corresponding to $u=(u_1(\eta t), u_2(\eta t), \dots)$, for some $\eta\in\Fq^*$, and assume that $C_\eta$ is real.

Let $e$ be the largest power of $2$ such that all of the $u_i$ are contained in $\Fq[t^\frac{e}2]$. Let $f$ be the order of $\eta$. If $|f|_2< e$ then $C_\eta=C$. If $|f|_2=e$ then $C_\eta$ is the conjugacy class corresponding to $v=(v_1(t), v_2(t), \dots)$, where $v_i(t)$ is the same as $u_i(t)$ except that the coefficient of $t^{ae}$ has reversed sign for odd $a$. If $|f|_2>e$ then $C_\eta$ is not real (resp. $\zeta$-real) which is a contradiction.

Thus there is at most one other real (resp. $\zeta$-real) conjugacy class in $GL_n(q)$ which is equivalent to $C$.
\end{proof}

\begin{lemma}\label{l: options}
let $C$ be a conjugacy class in $GL_n(q)$ of type $\nu$ with corresponding sequence $u=(u_1(t),u_2(t),\ldots)$. Let $d=|(n_1,n_2,\dots)|_2$.
\begin{enumerate}
\item If $d=1$, then the set of real conjugacy classes in $GL_n(q)$ is partitioned into equivalence classes of size $2$, and there are no $\zeta$-real classes.
\item Suppose that $d\geq 2$.
\begin{enumerate}
\item If $C$ is a real (resp. $\zeta$-real) class then $C$ is equivalent to one other real (resp. $\zeta$-real) class provided at least one $u_i$ is not in $\Fq[t^{|q-1|_2}]$. Moreover, if $C$ is a real class such that not all $u_i$ lie in $\Fq[t^{|q-1|_2}]$ then $C$ is not equivalent to any $\zeta$-real class.
\item If $C$ is a real (resp. $\zeta$-real) class, and all $u_i$ lie in $\Fq[t^{|q-1|_2}]$, then $C$ is equivalent to exactly one $\zeta$-real (resp. real) class; moreover $C$ is not equivalent to any other real (resp. $\zeta$-real) class.
\end{enumerate}
\end{enumerate}
\end{lemma}
\begin{proof}
Suppose first that $d=1$. Then $n_i$ is odd for some $i$ and, in particular, there are no $\zeta$-real classes of type $\nu$. Now suppose that $C$ is real. Then, since $n_i$ is odd, we have two distinct equivalent sequences which both correspond to a real class in $GL_n(q)$:
$$(u_1(t),u_2(t),\ldots) \textrm{ and } (u_1(-t),u_2(-t),\ldots).$$
Thus $C$ is equivalent to a distinct real class in $GL_n(q)$ as required.

Now suppose that $d\geq 2$, so that $n_i$ is even for all $i$. Let $C$ be a real (resp. $\zeta$-real) class such that at least one $u_i$ is not in $\Fq[t^{|q-1|_2}]$. Let $e$ be the largest power of $2$ such that all of the $u_i$ are contained in $\Fq[t^\frac{e}2]$ and take $\eta\in\Fq^*$ of order $e$. Then, once again, we have two distinct equivalent sequences which both correspond to a real (resp. $\zeta$-real) class in $GL_n(q)$:
$$(u_1(t),u_2(t),\ldots) \textrm{ and } (u_1(\eta t),u_2(\eta t),\ldots).$$
Thus $C$ is equivalent to one other real (resp. $\zeta$-real) class as required. Lemma \ref{l: lift2} implies that $C$ cannot be equivalent to a real (resp. $\zeta$-real) class.

Now suppose that $d\geq 2$ that all of the $u_i$ lie in $\Fq[t^{|q-1|_2}]$. If $C$ is real (resp. $\zeta$-real) then Lemma \ref{l: lift2} implies that $C$ is not equivalent to any other real (resp. $\zeta$-real) class. If $C$ is real (resp. $\zeta$-real) then define $\kappa$ (resp. $\theta$) just as in the proof of Lemma \ref{l: lift2}; the class corresponding to $\kappa u$ (resp. $\theta u$) is $\zeta$-real (resp. real).
\end{proof}

This lemma allows us to write down a formula for $pgl_\nu$ in all cases. Recall that $gl_\nu=\prod_{n_i>0} n_{q,n_i}$ and $n_{q,n_i}$ is defined in Lemma \ref{l: nqd}.

\begin{corollary}\label{c: pglnu}
Let $\nu=1^{n_1}2^{n_2}\cdots$ be a partition of $n$. If $d=1$ then $pgl_\nu = \frac12 gl_\nu$. If $d>1$ then $pgl_\nu=gl_\nu$.
\end{corollary}

\subsubsection{Conclusion}

We summarise our results for both odd and even characteristic in the following theorem:

\begin{theorem}\label{t: pglnq}
The number of real conjugacy classes in $PGL_n(q)$ is given by
$$\sum_{|\nu|=n}\frac{1}{2^{\sigma_{\nu}}}\prod_{n_i>0}n_{q,n_i}.$$
Here we set $d=|(n_1,n_2,\dots)|_2$ and define $\sigma_\nu$ to equal $0$ if $|dq|_2>1$ and to equal $1$ otherwise. Furthermore we note that all real conjugacy classes in $PGL_n(q)$ are strongly real.
\end{theorem}

We have not proved the statement about strong reality, however this follows from the work of Vinroot \cite[Theorem 3]{vinroot}. Note too that, for $q$ even or $n$ odd, the number of $GL_n(q)$-classes of $SL_n(q)$-real elements in $SL_n(q)$ is the same as the number of $PGL_n(q)$-real classes. This is reminiscent of an observation of Lehrer who pointed out that the total number of conjugacy classes in $PGL_n(q)$ is the same as the total number of $GL_n(q)$-classes in $SL_n(q)$ \cite{lehrer}.

\section{$PSL_n(q)$, $q$ is even or $|n|_2\neq |q-1|_2$}\label{s: psl1}

In this section we begin work on a classification of the real conjugacy classes in $PSL_n(q)$. We think of $PSL_n(q)$ as the image of $SL_n(q)$ under the quotient map $GL_n(q)\rightarrow PGL_n(q)=GL_n(q)/Z$ hence the elements of $PSL_n(q)$ will be written as $gZ$ where $g\in SL_n(q)$.

Recall that we may describe $gZ$ as being of type $\nu$, for some partition $\nu$, since all the elements in $GL_n(q)$ to which $gZ$ lifts are of the same type $\nu$. Our first result holds for any $q$ and for any cover of $PSL_n(q)$.

\begin{proposition}\label{p: pslsplit}
Let $Y\leq Z(GL_n(q))$ and set $G$ to equal $GL_n(q)/Y$ and $H$ to equal $SL_n(q)/(SL_n(q)\cap Y)$. Let $C$ be a $G$-conjugacy class in $H$ which is associated with a partition $\nu=(\nu_1,\nu_2,\dots)$. Then $C$ splits into $h_\nu=(q-1,\nu_1,\nu_2, \dots)$ $H$-conjugacy classes.
\end{proposition}
\begin{proof}
We use the methods of \cite{wall}. Let $C=C_1\cup\cdots \cup C_h$, where $C_1, \dots, C_h$ are $H$-conjugacy classes; to prove the proposition we must calculate the value of $h$. Suppose that $g$ projects into $C_1$ and let
$$X=\{h\in GL_n(q): hgh^{-1} = yg, \exists y\in Y\}.$$
Set $D_Y$ to be the group $\det X$ which lies in $\Fq^*.$

The class $C_1$ is stabilized in $G$ by $XH$. Now there is an isomorphism $G/(XH)\cong \Fq^*/D_Y$. Thus we have the equality $h=\frac{q-1}{|D_Y|}.$

We must now calculate $|D_Y|$. Suppose that $Y_1\leq Y_2 \leq Z(GL_n(q))$. It is clear that $D_{Y_1}\leq D_{Y_2}$ if and only if $Y_1\leq Y_2$. But now \cite[Theorem 4]{wall} states that $D_{Y_1}=D_{Y_2}$ for $Y_1=\{1\}$ and $Y=Z(GL_n(q))$. Hence $D_{Y_1}=D_{Y_2}$.

Now we have, in effect, already calculated $D_Y$ for $Y=\{1\}$; this is the content of \cite[(3.1)]{macdonald} and is described in the first paragraph of Section \ref{s: sl}. Thus we know that, for any $Y\leq Z(GL_n(q))$, $D_Y$ is the subgroup of $\Fq^*$ which has index $h_\nu=(q-1,\nu_1,\nu_2, \dots)$; the result follows.
\end{proof}

Now, similarly to before, we set $psl_\nu$ to equal the number of $PGL_n(q)$-real $PGL_n(q)$-conjugacy classes of type $\nu$ contained in $PSL_n(q)$.

Proposition \ref{p: pslsplit}, together with Corollary \ref{c: not2mod4} implies that, for $|n|_2\neq |q-1|_2$, the number of real conjugacy classes in $PSL_n(q)$ is equal to
$$\sum_{|\nu|=n} h_\nu psl_\nu.$$
Thus, for the remainder of this section (and the next), we will calculate $psl_\nu$ for differing $\nu, q$ and $n$.

\subsection{$q$ is even}

We know that all real elements in $PGL_n(q)$ lift to real elements in $GL_n(q)$. Since there are no equivalences for $q$ even, there is a 1-1 correspondence of real conjugacy classes between the two groups - indeed the same holds for $GL_n(q)/Y$ where $Y$ is any subgroup of $Z(GL_n(q))$.

Now all real elements in $GL_n(q)$ are in $SL_n(q)$. What is more these elements are real, in fact strongly real, in $SL_n(q)$. Hence there is also a 1-1 correspondence between real elements in $SL_n(q)$ and those in $PSL_n(q)$. We conclude that $psl_\nu = sl_\nu =  \prod_{n_i>0} n_{q,n_i}$.

\subsection{$q$ is odd and $|n|_2<|q-1|_2$}

\begin{lemma}\label{l: projreal2}
Suppose that $2\leq |n|_2<|q-1|_2$. If $g$ is $\zeta$-real in $GL_n(q)$ then $g$ does not project into $PSL_n(q)$.
\end{lemma}
\begin{proof}
Suppose that $hgh^{-1} = \zeta g$. Then $\det g = \pm \zeta^{\frac{n}2}$. Now take $\alpha\in\Fq^*$ and observe that $\det g(\alpha I) = \pm\zeta^{\frac{n}2}\alpha^n$. We may suppose without loss of generality that $\zeta$ generates $\Fq^*$ and suppose that $\alpha = \zeta^a$ for some integer $a$. Then $\det(\alpha g)= \zeta^{\frac{n}{2}+an}$ or $\zeta^{\frac{n}{2}+an+\frac{q-1}{2}}$. If $g$ projects down to $PSL_n(q)$ then we must have $\det(\alpha g)=1$ for some $\alpha\in\Fq^*$. But then we have $(q-1)|(\frac{n}{2}+an)$ or $(q-1)|(\frac{n}{2}+an+\frac{q-1}{2})$ which is impossible.
\end{proof}

\begin{lemma}
If $|n|_2<|q-1|_2$ then all real elements in $GL_n(q)$ project into $PSL_n(q)$. What is more these projections are strongly real in $PSL_n(q)$.
\end{lemma}
\begin{proof}
Lemma \ref{l: solxn} implies that there exists $\alpha\in\Fq$ such that $\alpha^n=-1$, thereby implying that $\det(\alpha I)=-1$. Now take $g$ real in $GL_n(q)$ so that, in particular, $\det g=\pm1$. If $\det g=-1$ then $\det(\alpha g)=1$ and $\alpha g$ is conjugate to $\alpha g^{-1}$. Thus $\alpha gZ$ lies in $PSL_n(q)$ as required.

Take $h\in GL_n(q)$ such that $hgh^{-1}=g^{-1}$. Corollary \ref{c: not2mod4} implies that all real elements in $GL_n(q)$ are strongly real in $\langle SL_n(q), (\alpha I) \rangle$, hence we may assume that $\det h = \pm 1$ and $h^2=1$. If $\det h=1$ then $hZ\in PSL_n(q)$ and so $gZ$ is strongly real in $PSL_n(q)$. If $\det h=-1$ then $\alpha hZ\in PSL_n(q)$ and, once more, $gZ$ is strongly real in $PSL_n(q)$.
\end{proof}

These two lemmas imply that the number of real classes in $PSL_n(q)$ is equal to the number of equivalence classes of real elements in $GL_n(q)$. We conclude that, in this case, $psl_\nu=\frac12 gl_\nu=\frac12 \prod_{n_i>0} n_{q,n_i}$.

\subsection{$q$ is odd and $|n|_2>|q-1|_2$.}\label{s: ngreater}

\begin{lemma}\label{l: det}
Suppose that $q$ is odd and $|n|_2\geq|q-1|_2$. If $g$ and $\eta g$ are both real (resp. both $\zeta$-real) then $\det g = \det(\eta g)$.
\end{lemma}
\begin{proof}
Observe that $\det(\eta g)= \eta^n \det g$. Now if $g$ and $\eta g$ are both real (or both $\zeta$-real) then $\det(\eta g) = \pm \det g$ and so $\eta^n=\pm 1$. Since $|n|_2\geq|q-1|_2$, Lemma \ref{l: solxn} implies that $\det h = \det g$.
\end{proof}

\begin{lemma}
Suppose that $q$ is odd and $|n|_2\geq|q-1|_2$. A $PGL_n(q)$-real conjugacy class $gZ\in PGL_n(q)$ is contained in $PSL_n(q)$ if and only if
\begin{enumerate}
\item it lifts to a real element in $GL_n(q)$ which is contained in $SL_n(q)$, or
\item it lifts to a $\zeta$-real element in $GL_n(q)$ of determinant $\zeta^{\frac{n}{2}}$ in the case where $|n|_2>|q-1|_2$, and of determinant $-\zeta^{\frac{n}{2}}$ otherwise.
\end{enumerate}
\end{lemma}
\begin{proof}
Suppose that $gZ$ lies in $PGL_n(q)$ with $g\in GL_n(q)$. If $g$ has determinant $-1$ in $GL_n(q)$ then $\det(\eta g) = \eta^n \det g = -\eta^n$, for $\eta\in \Fq^*$. Lemma \ref{l: solxn} implies that $\det(\eta g)\neq 1$ and so $g$ does not project into $PSL_n(q)$. Thus if $gZ$ is $PGL_n(q)$-real contained in $PSL_n(q)$ then there are two possibilities:
\begin{enumerate}
\item $g$ is real in $GL_n(q)$ with $\det g=1$;
\item $g$ is $\zeta$-real in $GL_n(q)$.
\end{enumerate}
In the second case, we can take $\zeta$ to be any non-square in $\Fq$; in particular we assume that $\zeta$ is a generator of the cyclic group $\Fq^*$. Now $g$ is conjugate to $\zeta g^{-1}$ and so $(\det g)^2 = \zeta^n$. In particular $\det g = \zeta^{\frac{n}2}$ or $\zeta^{\frac{n+q-1}2}$. Now $g$ projects into $PSL_n(q)$ only if there exists $\alpha\in\Fq^*$ such that $\alpha g\in SL_n(q)$. Write $\alpha=\zeta^b$ for some integer $b$. Then such an $\alpha$ exists provided one of the following equations has a solution:
$$\frac{n}2+bn\equiv 0\pmod{q-1}, \quad \frac{n+q-1}2+bn\equiv 0\pmod{q-1}.$$
These equations translate into two cases:
\begin{enumerate}
\item[(2a)] If $|n|_2>|q-1|_2$ then only the first solution is possible. This corresponds to the situation where $\det g= \zeta^{\frac{n}2}$.
\item[(2b)] If $|n|_2=|q-1|_2$ then only the second solution is possible. This corresponds to the situation where $\det g= -\zeta^{\frac{n}2}$.
\end{enumerate}
\end{proof}

For $\nu=1^{n_1}2^{n_2}\cdots$, set $d=|(n_1,n_2,\dots )|_2$ and, as before, define $\sigma_\nu$ to equal $0$ if $|dq|_2>1$ and to equal $1$ otherwise.

If $|n|_2>|q-1|_2$ and $d>1$ then the number of $\zeta$-real classes which project into $PSL_n(q)$ is the same as the number of real classes which project into $PSL_n(q)$; hence $psl_\nu=sl_\nu$. If $d=1$ then there are no $\zeta$-real classes in $GL_n(q)$ and so $psl_\nu = \frac12 sl_\nu$. Furthermore by Corollary \ref{c: not2mod4} we know that all of these conjugacy classes are strongly real in $PSL_n(q)$.

\subsection{Conclusion}

We summarise our findings in the following theorem:

\begin{theorem}\label{t: psl}
Suppose that $q$ is even or $|n|_2\neq |q-1|_2$. Let $d=|(n_1,n_2,\dots)|_2$. Then the number of real classes in $PSL_n(q)$ of type $\nu$ is equal to $h_\nu psl_\nu$ where
$$psl_\nu= \left\{ \begin{array}{ll}
\frac{1}{(2,q-1)} \prod_{n_i>0} n_{q,n_i}, & |n|_2<|q-1|_2 \textrm{ or } q \textrm{ is even}; \\
sl_\nu, & |n|_2>|q-1|_2, d>1 \textrm{ and } q \textrm{ is odd}.\\
\frac12 sl_\nu, & |n|_2>|q-1|_2, d=1 \textrm{ and } q \textrm{ is odd}. \\
          \end{array}\right.
$$
What is more all real classes in $PSL_n(q)$ are strongly real.
\end{theorem}

Note that, when $q$ is even, Theorem \ref{t: psl} also holds for $SL_n(q)/(SL_n(q)\cap Y)$ where $Y$ is any subgroup of $Z(GL_n(q))$.

\section{$PSL_n(q)$, $q$ is odd and $|n|_2=|q-1|_2$}\label{s: psl2}

Now as before we set $psl_\nu$ to equal the number of $PGL_n(q)$-real $PGL_n(q)$-conjugacy classes of type $\nu$ contained in $PSL_n(q)$. We start by calculating $psl_\nu$ for various scenarios. Note that both of the lemmas from Section \ref{s: ngreater} apply here. 

First set $\nu=(1^{n_1}2^{n_2}\ldots)$, $d=|(n_1,n_2,\dots)|_2$. We will use the methods and notation of Proposition~\ref{p: slnu}, and consider various cases. In particular suppose that $C$ is a real class (resp. a $\zeta$-real class) of type $\nu$ in $GL_n(q)$. Then $C$ is associated with a sequence of polynomials $(u_1(t), u_2(t), \dots)$ such that the $u_i(t)$ are self-reciprocal (resp. $\zeta$-self-reciprocal). Now let $a_i$ be the leading term in $u_i(t)$; then, for $g\in C$, $\det g = (-1)^n\prod_{n_i>0} a_i^i= \prod_{n_i>0} a_i^i$. We know that $a_i$ is equal to $\pm 1$ (resp. $\pm\frac1{\zeta^{\frac{n_i}2}})$ when $C$ is real (resp. $C$ is $\zeta$-real).

For this section it will help to choose $\zeta$ to be a non-square which satisfies $\zeta^{\frac{n}2}=-1$; this means that a $\zeta$-real element, like a real element, will have determinant $\pm 1$. Then, for $g$ real or $\zeta$-real, Lemma \ref{l: solxn} implies that, if $\det g=-1$, then $g$ does not project into $PSL_n(q)$.

\begin{enumerate}
\item[(P1)]\label{p1} Suppose that $d=1$; thus there are no $\zeta$-real elements. We must have $\det g=1$ and so Lemma \ref{l: options} implies that $C$ is equivalent to one other real class in $GL_n(q)$; Lemma \ref{l: det} implies that this class consists of elements of determinant $1$. Thus $psl_\nu=\frac12 sl_\nu$. 

\item[(P2)]\label{p2} Suppose that $d>1$ and that $n_i=0$ for all odd $i$; in particular this means that $d<|n|_2$. If $C$ is real then Lemma \ref{l: options} implies that $C$ is equivalent to one other real class in $GL_n(q)$; Lemma \ref{l: det} implies that this class consists of elements of determinant $1$.

If $C$ is $\zeta$-real then we must have $d>1$ and 
$$\det g= \prod_{n_i>0}a_i^i = \prod_{n_i>0}\left(\pm\frac1{\zeta^{\frac{n_i}2}}\right)^i=  \frac1{\zeta^{\frac{n}2}}.$$
Thanks to our choice of $\zeta$ we have $\det g= \zeta^{\frac{n}2}=-1$ and so, as we have already observed, $g$ does not project into $PSL_n(q)$. Hence once again we have $psl_\nu=\frac12 sl_\nu.$

\item[(P3)]\label{p3} Suppose that $d>1$ and that there exists $i$ odd for which $n_i>0$. The number of real classes which lie in $SL_n(q)$ is given by Proposition \ref{p: slnu} and is equal to
$$f_\nu(q)\prod_{i \textrm { odd}, n_i>0} q^{\frac{n_i}2-1}\times \prod_{i \textrm { even},n_i>0}n_{q,n_i}.$$

We also need to count the number of $\zeta$-real classes for which the determinant is equal to $-\frac1{\zeta^{\frac{n}2}}=1$. The same methods as in Proposition \ref{p: slnu} yield that the number of such classes is equal to
$$g_\nu(q)\prod_{i \textrm { odd}, n_i>0} q^{\frac{n_i}2-1}\times \prod_{i \textrm { even},n_i>0}n_{q,n_i}$$
where $g_\nu(q)=\frac{(q+1)^r-(q-1)^r}{2}$ where $r$ is the number of odd values of $i$ for which $n_i>0$.

The total number of these two types of conjugacy class is $\prod_{n_i>0}n_{q,n_i}$. Lemma \ref{l: options} implies that these conjugacy classes partition into $psl_\nu=\frac12\prod_{n_i>0}n_{q,n_i}$ equivalence classes.
\end{enumerate}

Note that if $4|n$ then Proposition \ref{p: not2mod4} implies that all of the above classes are strongly real in $PSL_n(q)$. Hence, using Proposition \ref{p: pslsplit}, we have the following:

\begin{proposition}
Suppose that $|n|_2=|q-1|_2$ and $4|n$. Then the total number of real conjugacy classes in $PSL_n(q)$ is the same as the total number of strongly real conjugacy classes and is given by
$$\sum_{|\nu|=n}h_\nu psl_{\nu}$$
where the values for $psl_\nu$ are as outlined above.
\end{proposition}

\subsection{$n\equiv 2\pmod 4$ and $q\equiv 3\pmod 4$}

This is the only case left to consider for $PSL_n(q)$. In the five points above we have calculated the number of $PGL_n(q)$-classes of $PGL_n(q)$-classes lying in $PSL_n(q)$. But in this case we do not know if all of these classes will remain real in $PSL_n(q)$.

\begin{proposition}
Suppose that $n\equiv 2\pmod 4$ and $q\equiv 3\pmod 4$. Then the total number of real conjugacy classes in $PSL_n(q)$ is given by
$$\sum_{|\nu|=n}h_\nu psl_{\nu}$$
where $psl_\nu$ is non-zero exactly when $n_i>0$ for some odd $i$. In this case the values for $psl_\nu$ are given by (P1), (P2) and (P3).
\end{proposition}
\begin{proof}
If $n_i>0$ for some odd $i$, then Proposition \ref{p: n2mod4} implies that a real (or a $\zeta$-real) conjugacy class of $GL_n(q)$ contained in $SL_n(q)$ is real (or $\zeta$-real) within $SL_n(q)$.

Hence we only need to deal with the situations of (P1) and (P2) where $n_i=0$ for all odd $i$. In fact (P2) cannot occur for $n\equiv 2\pmod 4$.

Thus we are left with the case (P1) only and there are no $\zeta$-real elements. Furthermore, for a real element, Lemma \ref{l: square} and Proposition \ref{p: n2mod4} imply that all reversing elements have non-square determinant in $\Fq$. But such elements do not project into $PSL_n(q)$ (see Lemma \ref{l: solxn}), hence this situation does not yield real elements in $PSL_n(q)$.
\end{proof}

\subsection{Conclusion}

We summarise our results in the following theorem.

\begin{theorem}\label{t: pslconc}
Suppose that $|n|_2=|q-1|_2$ and let $d=|(n_1,n_2,\dots)|_2$. Then the number of real classes in $PSL_n(q)$ is given by
$$\sum_{\nu}h_\nu psl_\nu.$$
If $4|n$ then,
$$psl_\nu=\left\{\begin{array}{ll}
\frac12 \prod_{n_i>0} n_{q,n_i}, & \textrm{if }d>1 \textrm{ and there exists }i \textrm{ odd for which } n_i>0 \\
\frac12sl_\nu, & \textrm{otherwise.}
\end{array}
 \right.$$
If $4\not |n$ then,
$$psl_\nu=\left\{\begin{array}{ll}
\frac12 \prod_{n_i>0} n_{q,n_i}, & \textrm{if }d>1 \textrm{ and there exists }i \textrm{ odd for which } n_i>0 \\
\frac12sl_\nu, & \textrm{if }d=1 \textrm{ and there exists }i \textrm{ odd for which } n_i>0 \\
0, & \textrm{otherwise.}
\end{array}
\right.$$
\end{theorem}

\section{Strongly real classes in $PSL_n(q)$}\label{pslstrongly}

Provided $n\not\equiv 2\pmod 4$ or $q\not\equiv 3\pmod 4$ we have shown that reality and strong reality coincide in $PSL_n(q)$. Throughout this section we examine the strongly real classes in $PSL_n(q)$ when $n\equiv 2\pmod 4$ and $q\equiv 3\pmod 4$. We start with a lemma.

\begin{lemma}\label{l: prev}
Suppose that $n\equiv 2 \pmod 4$ and $q\equiv 3 \pmod 4$. An element $gZ$ is strongly real in $PSL_n(q)$ if and only if $gZ$ lifts to an element $g$ in $GL_n(q)$ for which there is an element $h$ satisfying
\begin{enumerate}
\item $hgh^{-1}=g^{-1}$ (or $\zeta g^{-1}$);
\item $h^2\in Z(GL_n(q))$;
\item $\det h$ is a square.
\end{enumerate}
\end{lemma}
\begin{proof}
Suppose that such an element $h$ exists. Then $h^c$ has determinant $1$ for some odd integer $c$. Furthermore $(h^c)^2\in Z(GL_n(q))$ and $(h^c)g(h^c)^{-1}=g^{-1}$ (or $\zeta g^{-1}$). Thus $gZ$ is strongly real in $PSL_n(q)$.

On the other hand if $gZ$ is strongly real in $PSL_n(q)$ then, by definition, an element $h$ exists in $GL_n(q)$ satisfying the first two criteria given. What is more $h$ projects into $PSL_n(q)$; in other words $\eta h$ has determinant $1$ for some scalar $\eta$. This means that $\det h = \frac{1}{\eta^n}$ which is a square since $n$ is even.
\end{proof}

Take $gZ$ real in $PSL_n(q)$. Then $g$ is of type $\nu$ where $n_i>0$ for some odd $i$. Let $g$ be a real or $\zeta$-real element in $GL_n(q)$ and let $V$ be the module associated with $g$. Let $h$ be a reversing element for $g$ in $GL_n(q)$ which satisfies $h^2\in Z(GL_n(q))$.

Now $h$ permutes the minimal cyclic submodules of $V$ with orbits of size $2$ (in the proof of Proposition \ref{p: n2mod4} we called these orbits {\it $h$-minimal} submodules of $V$). This fact allows us to break the general situation into smaller subcases which we deal with in the next two lemmas.

\begin{lemma}\label{l: prev1}
Suppose that $V= W_p\oplus W_q$ where $W_p$ and $W_q$ are cyclic modules with irreducible characteristic polynomials $p(t)^a$ and $\tilde p(t)^a$ (resp. $\breve p(t)^a$). Furthermore assume that $h$ swaps $W_p$ and $W_q$. Set the degree of $p(t)$ (and $q(t)$) to be $d$. Then $\det h$ can be a square or a non-square if $ad$ is odd; otherwise $\det h$ is a square.
\end{lemma}
\begin{proof}
We proceed as per the proof of Lemma \ref{l: sltwocycles}; in particular, we can take $g$ to equal $\left(\begin{array}{cc}
B & 0 \\
0 & B^{-1}  \end{array}\right)$ or $\left(\begin{array}{cc}
B & 0 \\
0 & \zeta B^{-1}  \end{array}\right)$ for some $B\in GL_{ad}(q)$. This means that $h=\left(\begin{array}{cc}
0 & X \\
Y & 0  \end{array}\right)$ where $X$ and $Y$ centralize $B$ in $GL_{ad}(q)$. Then
$$\det h=(-1)^{ad}(\det X)(\det Y)=(-1)^{ad}(\det XY).$$
Since $h^2\in Z(GL_n(q))$ we must have $XY\in Z(GL_{ad}(q))$  and so 
$$\det h=(-1)^{ad}\alpha^{ad}=(-\alpha)^{ad}$$ 
where $\alpha\in\Fq$. Thus if $ad$ is even this determinant is a square. On the other hand if $ad$ is odd then we can let $X=Y=I$ and $\det h$ is a non-square, or take $X=I=-Y$ and $\det h$ is a square.
\end{proof}

\begin{lemma}\label{l: prev2}
Suppose that $V=W_p$, a cyclic module $W_p$ with irreducible characteristic polynomial $p(t)^a;$ furthermore $p(t)$ is self-reciprocal (resp. $\zeta$-self-reciprocal). Set the degree of $p(t)$ to be $d$. 
\begin{enumerate}
\item If $d$ is odd then $\det h$ can be chosen to be a square or a non-square if $a$ is odd; otherwise $\det h$ is a square for $a\equiv 0\pmod 4$ and $\det h$ is a non-square for $a\equiv 2\pmod 4$.
\item If $d$ is even and $a$ is even, then $\det h$ is a square.
\item If $d$ is even and $a$ is odd, then $\det h$ can be chosen to be a square or a non-square if $d\equiv 2\pmod 4$; otherwise $\det h$ is a square for $d\equiv 0\pmod 4$.
\end{enumerate}
\end{lemma}
\begin{proof}
Let us examine the relevant cases. We will use Lemmas \ref{l: cyclic} and \ref{l: cyclic2}; these give conditions for $\det h$ to be a square, but they do not assume that $h^2\in Z(GL_n(q))$. In the cases where these lemmas allow for $\det h$ to be a square or a non-square we need to check the situation under this extra assumption.

Suppose that $d$ is odd. Then $p(t)=t\pm 1$ and we refer to Lemma \ref{l: cyclic} and observe that the conclusions given there apply here also. The only thing we have to check is that the element $h$ satisfies $h^2\in Z(GL_n(q))$. But the $h$ which we exhibit in the proof is an involution so we are done.

Suppose that $d$ is even and $a$ is even. Then Lemma \ref{l: cyclic2} implies that $\det h$ is a square.

Suppose that $d$ is even and $a$ is odd. Write the Jordan decomposition in $GL_{ad}(q)$: $g=g_sg_u$. Then $g_s$ is centralized by $GL_a(q^d)$ and so the centralizer of $g$ must lie in $GL_a(q^d)$. Furthermore a reversing element for $g$ must be a reversing element for $g_s$ and hence must normalize $C_G(g_s)$. Thus this element must act as a field involution of $GL_a(q^d)$. 

Suppose that $d\equiv 2\pmod 4$. Corollary \ref{c: not2mod4} implies that there exists a reversing involution $h_0$. Since $a$ is odd, we can choose $z\in Z(GL_a(q^d)) = Z(C_G(g_s))$ such that $\det z$ is a non-square. Now $h_0$ acts as a field automorphism on $C_G(g_s)$ hence $$(zh_0)^2 = zz^{h_0} h_0^2 = z^{q^{\frac{d}{2}}+1}.$$ 
Clearly $h_0$ and $zh_0$ are reversing elements with different determinant. Now write $(zh_0)^2$ as an element of $GL_a(q^d)$: $(zh_0)^2 = \beta I$ for some $\beta \in \Fqd$. For this to lie in $Z(GL_{ad}(q))$ we must have $\beta^{(q^{\frac{d}2}+1)(q-1)}=1$. 

Since $|q^{\frac{d}{2}}-1|_2 = |q-1|_2$ we can take an odd power of $z$, $z^c$ say, such that $(z^c h_0)^2\in Z(GL_{ad}(q))$. Clearly $z^c h_0$ is a reversing element for $g$ and $\det (z^ch_0)$ is a square if and only if $\det (zh_0)$ is a square. We conclude that, in this situation, we can take $h$ to have determinant a square or a non-square.

Suppose that $d\equiv 0\pmod 4$. By Corollary \ref{c: not2mod4} we know that a reversing involution $h_0$ exists which acts as a field automorphism on $GL_a(q^d)$ and has determinant a square. Then any other reversing element must have form $zh_0$ where $z$ centralizes $g$. The form of $z$ is given (as an element in $GL_a(q^d)$) by 
$$z=\left(\begin{array}{cccc}
\beta_1 & \beta_2&&  \\
&\ddots &\ddots& \\
&&\ddots & \beta_2 \\
&&& \beta_1
  \end{array}\right).$$
Then
$$(zh_0)^2=zz^{h_0} h_0^2 = zz^{h_0} = \left(\begin{array}{ccc}
\beta_1^{q^{\frac{d}2}+1} & &  \\
&\ddots& \\
&& \beta_1^{q^{\frac{d}2}+1}
  \end{array}\right).$$
For $(zh_0)^2$ to lie in $Z(GL_{ad}(q))$, we must have $\beta_1^{(q^{\frac{d}2}+1)(q-1)}=1$. But this means that $\beta_1$ must be a square in $\Fqd$. Thus $\det zh_0$ is a square in all cases.
\end{proof}

It is now just a matter of summing up what we have proved so far, and converting our result into the language of Macdonald.

\begin{theorem}\label{t: pslstrong}
Let $n\equiv 2\pmod 4$ and $q\equiv 3\pmod 4$. Let $gZ$ be real in $PSL_n(q)$ of type $\nu=1^{n_1}2^{n_2}\cdots$, and suppose that $g$ can be taken to be self-reciprocal (resp. $\zeta$-self-reciprocal). Then $gZ$ is {\bf not} strongly real in $PSL_n(q)$ if and only if the following conditions hold for all odd $i$ such that $n_i>0$:
\begin{enumerate}
\item All factors of $u_i(t)$ have even degree.
\item All self-reciprocal (resp. $\zeta$-self-reciprocal) factors of $u_i(t)$ have degree equivalent to $0\pmod 4$.
\end{enumerate}
\end{theorem}

Note that Theorem \ref{t: pslconc} implies that an odd $i$ exists for which $n_i>0$.

\begin{proof}
Let $V$ be the module associated with $g$. By Lemma \ref{l: prev} we need to show that any reversing element $h$, which satisfies $h^2\in Z(GL_n(q))$, has $\det h$ a non-square.

Break $V$ up into $h$-minimal submodules, $W$, in a similar way to the proof of Proposition \ref{p: n2mod4}. Suppose that $g|_W$ is reversible by an element $h_0$ which satisfies $h^2\in Z(GL(W))$ and for which we can choose $\det h_0$ to be square or non-square. Then clearly we can choose $h$ to be a square or a non-square; Lemmas \ref{l: prev1} and \ref{l: prev2} give the conditions under which this is possible. These conditions are precisely the ones excluded by the statement of the theorem.

Thus the conditions given in the theorem ensure that, for every $h$-minimal submodule $W$, the reversing elements of $g|_W$ in $GL(W)$ either all have determinant a square or all have determinant a non-square; in fact Lemmas \ref{l: prev1} and \ref{l: prev2} imply that the determinant will be a square if and only if the dimension of $W$ is divisible by $4$. Since $n\equiv 2\pmod 4$ we conclude that $\det h$ must have determinant a non-square as required.
\end{proof}

\section{Quotients of $SL_n(q)$}\label{s: isogenous}

We examine the real and strongly-real classes in $SL_n(q)/Y$ where $Y$ is some subgroup of $Z(SL_n(q))$. We have noted already that Theorem \ref{t: psl} holds for $SL_n(q)/Y$ where $q$ is even. In fact, if $|Y|$ is odd, then the number of real (resp. strongly real) classes will equal the number of real (resp. strongly real) classes in $SL_n(q)$. Similarly if $|Y|_2=|(n,q-1)|_2$ then the number of such classes will be the same as in $PSL_n(q)$. Hence, in this section, we assume that $1<|Y|_2<|(n,q-1)|_2$; in particular we assume that 
$$q\equiv1\pmod 4, \textrm{ and } n \equiv 0\pmod 4.$$

In what follows we will think of $Y$ as being a subgroup of $SL_n(q)$, $GL_n(q)$ or $\Fq^*$ depending on the context. We need two new concepts that mirror our treatment of projective groups from Section \ref{s: pgl}.

Firstly we define two elements, $g_1$ and $g_2$, of $GL_n(q)$ to be {\it $Y$-equivalent} if they project onto the same element of $GL_n(q)/Y$; so $g_2= g_1y$ for some $y\in Y$. This notion can be extended to conjugacy classes of $GL_n(q)$ and $GL_n(q)/Y$.

Secondly we generalize the idea of a $\zeta$-real element. Let $\zeta_Y$ be an element of $Y$ such that $\zeta_Y\neq \alpha^2$ for all $\alpha\in Y$ ($\zeta_Y$ is a {\it non-square} in $Y$); we say that $g$ is {\it $\zeta_Y$-real} in $GL_n(q)$ if there exists $h\in GL_n(q)$ such that $hgh^{-1} = \zeta_Y g^{-1}$. It is easy to see that all real elements in $SL_n(q)/Y$ will lift to a real element or a $\zeta_Y$-real element in $GL_n(q)$ (c.f. Lemma \ref{l: projreal}).

For ease of calculation we will set $\zeta_Y$ to be an element which satisfies $\zeta_Y^{|Y|_2} = 1$. In particular this means that all $\zeta_Y$-real elements, like all real elements, have determinant $\pm 1$. Since $|Y|_2<|n|_2,$ we know that only elements of determinant $1$ project into $SL_n(q)/Y$ (c.f. Lemma \ref{l: solxn}).

Now Proposition \ref{p: not2mod4} states that all $GL_n(q)$-real elements in $SL_n(q)$ are strongly real in $SL_n(q)$. It is easy enough to modify the proof to show that all $GL_n(q)$-$\zeta_Y$-real elements in $SL_n(q)$ are strongly $\zeta_Y$-real in $SL_n(q)$ (where {\it strongly $\zeta_Y$-real} has the obvious definition).

Finally Proposition \ref{p: pslsplit} implies that if a $GL_n(q)/Y$-class is of type $\nu$ then the class will split into $h_\nu$ classes in $SL_n(q)/Y$. All this combines to give the following proposition:

\begin{proposition}
Let $Y$ be a subset of $Z(SL_n(q))$ such that $1<|Y|_2<|(n,q-1)|_2$. The total number of real classes in $SL_n(q)/Y$ is the same as the number of strongly real classes in $SL_n(q)/Y$ and is equal to
$$\sum_{|\nu|=n}h_\nu sly_{\nu}. $$
Here $sly_\nu$ is the number of $Y$-equivalence classes in the set of all real and $\zeta_Y$-real conjugacy classes of type $\nu$ and determinant $1$ in $GL_n(q)$. 
\end{proposition}

All that remains is to calculate the value of $sly_\nu$ for differing $\nu$, $Y$, $q$ and $n$. In fact, though, this is easy. The number of real classes and the number of $\zeta_Y$-real conjugacy classes will be the same (c.f. Lemma \ref{l: nqdzeta} and note that $\zeta_Y$-self-reciprocal polynomials exist with odd degree).  These will be partitioned into sets of size $2$ as described in Lemma \ref{l: options}. Hence $sly_\nu = sl_\nu$ and the value of $sl_\nu$ is given in Proposition \ref{p: slnu}.

We summarise our results as follows.

\begin{theorem}
Let $Y$ be a subset of $Z(SL_n(q))$ such that $1<|Y|_2<|(n,q-1)|_2$. The number of real classes in $SL_n(q)/Y$ is equal to the number of strongly real classes, and is given by
$$\sum_{|\nu|=n}h_\nu sl_\nu .$$
\end{theorem}

Observe that this is the same as the number of real classes in $SL_n(q)$.

\section{Some exceptional cases}\label{s: exceptional}

In order to complete our classification of real classes in all quasi-simple covers of $PSL_n(q)$ we must deal with some exceptional situations, namely quasi-simple covers of $PSL_n(q)$ which are not quotients of $SL_n(q)$. There are five situations where this may occur:  $PSL_2(4)$, $PSL_3(2)$, $PSL_2(9)$, $PSL_3(4)$ and $PSL_4(2)$ \cite[Theorem 5.1.4]{kl}.

Write $M(G)$ for the Schur multiplier of a simple group $G$. If $G=PSL_2(4)$ (resp. $PSL_3(2)$) then $|M(G)|=2$ and the double cover of $G$ is isomorphic to $SL_2(5)$ (resp. $SL_2(7)$). We have already analysed the real classes in these groups. The remaining three groups need to be analysed in turn; we start by recording some information about each (see \cite[Proposition 2.9.1 and Theorem 5.1.4]{kl}):
\begin{center}
\begin{tabular}{|c|c|c|}
\hline
$G$ & Isomorphism & $M(G)$ \\
\hline
$PSL_2(9)$ & $A_6$ & $C_6$ \\
$PSL_4(2)$ & $A_8$ & $C_2$ \\
$PSL_3(4)$ & & $C_4\times C_{12}$ \\
\hline
\end{tabular}
\end{center}
Here $C_n$ is the cyclic group of order $n$, and the middle column lists groups to which $G$ is isomorphic. Information about conjugacy and reality can, for quasi-simple groups with cyclic centre, be found in \cite{atlas}; however classifying the strongly real classes is more tricky, so we prefer to calculate our results directly, where possible.

Note that our approach in this section is, in some sense, the reverse of what has happened in the rest of the paper. In this section we have complete information about (strongly) real classes in our simple group $G$, and we wish to deduce information about (strongly) real classes in quasi-simple covers of $G$. We start with a lemma which applies to this situation in some generality.

\begin{lemma}\label{l: lifting}
Let $G, H$ be groups such that $H/Z\cong G$ where $Z$ is an odd-order central subgroup of $H$. Let $C$ be a real class in $G$ containing elements of order $n$; then $C$ lifts to a unique real class $C_H$ in $H$ and this class consists of elements of order $n$. What is more if $C$ is strongly real than $C_H$ is strongly real.
\end{lemma}
\begin{proof}
Let $\chi$ be a real-valued irreducible complex character (or {\it rvicc}) of $H$. Let $g\in Z$; then $\chi(g)=1$ for $g\neq 1$. This implies that $\chi$ is an rvicc of $G$. Since every rvicc of $G$ is an rvicc of $H$, we conclude that $G$ and $H$ have the same number of rvicc's. Thus $G$ and $H$ have the same number of real classes.

Now suppose that $h$ and $gh$ lie in different conjugacy classes of $H$, with $g\in Z, g\neq 1$. Let $\rho$ be the regular character; then $\rho(g)=0$ \cite[(2.10)]{isaacs}, hence there exists an irreducible complex representation $\phi$, of dimension $n$, such that $\phi(g)\neq I$. Then, since the order of $g$ is odd, Schur's Lemma implies that $\phi(g)=\eta I$ for some $\eta\not\in\mathbb{R}$. Let $\phi_1$ be the character ${\rm Tr}(\phi)$; then $\phi_1(gh)=n\eta\phi_1(h)$ and so $h$ and $gh$ cannot both be real; thus every real class of $G$ lifts to a unique real class in $H$.

Suppose next that $h$ is real in $H$ and $hZ$ is real in $G$ of order $m$. Then $h^m\in Z$, and the only element in $Z$ that is real in $H$ is the identity. Thus $h^m=1$ as required.

Finally suppose that $C_H$ is a real class in $H$ such that $h\in C_H$ and $hZ$ lies in $C$, a strongly real class in $G$. Then there exists $fZ\in G$ such that
$$(fZ)(hZ)(fZ)^{-1}=h^{-1}Z, {\rm \ and \ } (fZ)^2 = Z\in G.$$
Since $Z$ has odd order we can assume that $f^2=1\in H$. Then $fhf^{-1} = h^{-1}g$ for some $g\in Z$. Since $Z$ has odd order, $g=g_1^{-2}$ for some $g_1\in Z$. Then $f(hg_1)f^{-1} = (hg_1)^{-1}.$ Thus $hg^{-1}$ is strongly real in $H$ and projects onto $hZ$ in $G$. Since $C_H$ is the unique real class to which $C$ lifts, we conclude that $C_H$ is strongly real.
\end{proof}

Note that a consequence of this lemma is that $G$ and $H$ have the same number of real (resp. strongly real) classes. If $G$ is isomorphic to either $PSL_2(9)$ or  $PSL_2(4)$, then $M(G)$ contains a unique involution; thus $G$ has a unique double cover, $2.G$. Lemma \ref{l: lifting} reduces the problem of studying real and strongly real classes in the covers of $G$ to the problem of studying the real and strongly real classes in $G$ and $2.G$. The case of $PSL_3(4)$ is more difficult.

Throughout what follows $G$ is a simple group, $H$ a quasi-simple group with centre $Z=Z(H)$ such that $G\cong H/Z$.

\subsection{Covers of $PSL_2(9)$}

The group $PSL_2(9)\cong A_6$ contains 7 conjugacy classes (with elements of order 1,2,3,3,4,5, and 5), all of which are real (see Theorem \ref{t: psl} or \cite{tz}). Lemma \ref{l: lifting} implies that $3.PSL_2(9)$ has 7 real conjugacy classes with elements of the same orders. Theorem \ref{t: psl} implies that all of the conjugacy classes in $PSL_2(9)$ are strongly real, hence the same is true of $3.PSL_2(9)$.

Similarly $SL_2(9)\cong 2.PSL_2(9)$ contains 13 conjugacy classes (with elements of orders 1,2,3,3,4,5,5,6,6,8,8,10, and 10) all of which are real. Lemma \ref{l: lifting} implies that $6.PSL_2(9)$ has 13 real conjugacy classes with elements of the same orders.

The groups $SL_2(9)\cong 2.PSL_2(9)$ and $6.PSL_2(9)$ contain a single involution (in the centre); thus both groups contain precisely two strongly real classes.

\subsection{Covers of $PSL_4(2)$}
Let $H$ be a double cover of $A_{2n}$. Consider the group $H_{2n}$ given by the presentation
$$H_{2n}=\langle s_1, \dots s_{2n-1} \mid s_k^2, -(s_ks_{k+1})^3, -(s_ks_j)^2\rangle,$$
where $j,k=1,\dots, 2n-1$, $|j-k|>1$, and $-1$ is defined to be a central element. The group $H_{2n}$ is a double cover of $S_{2n}$ with center $Z=\{1, -1\}$ \cite[p.175]{aschbacher3}; the projection map is given by 
$$\pi: H_{2n}\mapsto S_{2n}, s_k \mapsto (k \, \, {k\!\!+\!\!1}).$$ 
Then $H_{2n}$ has a subgroup $J_{2n}$ of index $2$ which is the double cover of $A_{2n}$; clearly $J_{2n}$ consists of all elements $\pm x$ where $x$ is the product of an even number of the $s_i$.

Now let $g$ be a real element of $J_{2n}$; then $gZ$ is a real element of $A_{2n}$. What is more if $g$ is (strongly) real in $J_{2n}$ then $-g$ is also (strongly) real (since $hgh^{-1} = g^{-1}$ implies that $h(-g)h^{-1} = -g^{-1} = (-g)^{-1}$).

With this in mind we list the real and strongly real classes in $J_8\cong 2.PSL_4(2)$ in the following table.

\begin{center}
\begin{tabular}{|c|c|c|c|c|c|}
\hline
Line & $gZ\in A_8$ & $g\in J_8$ & Order in $J_8$ & Real & Strongly Real\\
\hline
1 & $(1)$ & $\pm 1$ & $1, 2$ & Yes & Yes\\
2 & $(12)(34)$ & $\pm s_1s_3$ & $4$ & Yes & No\\
3 & $(123)$ & $\pm s_1s_2$ & $3,6$ & Yes & No\\
4 & $(123)(456)$ & $\pm s_1s_2s_4s_5$ & $3 , 6$ & Yes & Yes\\
5 & $(1234)(56)$ & $\pm s_1s_2s_3s_5$ & $8$ & Yes & No\\
6 & $(12345)$ & $\pm s_1s_2s_3s_4$ & $5,10$ & Yes & No\\
7 & $(12)(34)(56)(78)$ & $\pm s_1s_3s_5s_7$ & $2$  & Yes & Yes\\
8 & $(1234)(5678)$ & $\pm s_1s_2s_3s_5s_6s_7$ & $4$ & Yes & Yes\\
9 & $(123)(45)(67)$ & $\pm s_1s_2s_4s_6$ & $12$ & Yes & No\\
10 & $(123456)(78)$ & $\pm s_1s_2s_3s_4s_5s_7$ & $6,6$ & No & No\\
\hline
\end{tabular}
\end{center}

We need to explain the columns of this table: The first column records the line number. The second column lists representatives from all of the real classes in $A_8$. The third column lists the two elements in $J_8$ that project onto the given representative in $A_8$. The fourth column gives the order of elements in $J_8$ which project onto $gZ$ in $A_8$; the presence of two numbers in this column means that there are two different conjugacy classes of elements in $J_8$ that project onto the same conjugacy class in $A_8$. The final two columns state whether or not the elements $\pm g$ are (strongly) real.

\begin{proposition}
The above table is correct.
\end{proposition}
\begin{proof}
It is easy to check that the listed elements $gZ$ really do represent all of the real conjugacy classes in $A_8$ (one could check this directly, or else use the isomorphism with $PSL_4(2)$ and our theory above). That the elements $g$ project onto $gZ$ is also easy to check.

Now let us explain the third column. When $gZ$ has odd order $k$, it is clear that the set $gZ\subset 2.A_8$ must contain an element of order $k$ and an element of order $2k$, hence lines 1,3,4 and 6 are justified. (Indeed it is easy to establish the order of the elements $g$ directly.)

When $gZ$ has even order it is an easy matter to establish the order of $g$; we need only establish if there are one or two conjugacy classes in each case. Observe that
\begin{equation*}
\begin{aligned}
(s_7s_5)(s_1s_2s_3s_5)(s_7s_5)^{-1} &= -s_1s_2s_3s_5; \\
(s_1s_2s_1s_2s_3s_2)(s_1s_3s_5s_7)(s_1s_2s_1s_2s_3s_2)^{-1} &=-s_1s_3s_5s_7; \\
(s_5s_6s_5s_4s_5s_4)(s_1s_2s_4s_6)(s_5s_6s_5s_4s_5s_4)^{-1} &= -s_1s_2s_4s_6.
\end{aligned}
\end{equation*}

These justify lines 5, 7 and 9 respectively; line 2 is a consequence (consider $g^3$ for $g$ in line 9). Line 8 follows from the fact that
$hs_1s_2s_3s_5s_6s_7h^{-1} = -s_1s_2s_3s_5s_6s_7$,
where $h$ is any element which lifts to $(15)(26)(37)(48)$ in $A_8$.

Finally, for the last line, observe that $C_{A_8}(gZ)=\langle gZ \rangle$. Any element $h\in J_8$ which maps $g$ to $-g$ must project onto an element of $C_{A_8}(gZ)=\langle gZ \rangle$. All such elements $h$ commute with $g$, hence there are two conjugacy classes of this form.

Now let us consider the final two columns. Recall that, since $g$ is (strongly) real if and only if $-g$ is (strongly) real, we need only prove the result for any $h$ projecting onto $gZ$.

We start with the fourth column. It is clear that if a conjugacy class $C$ is real in $A_8$ and it lifts to a single conjugacy class in $J_8$ then this class is real in $J_8$. Similarly if $C$ lifts to two conjugacy classes containing elements of different order, then both of these classes are real in $J_8$. Thus the only classes left to consider are those on line 10.

Now observe that
$$(s_7s_4s_3s_4s_5s_4s_3s_1)(s_1s_2s_3s_4s_5s_7)(s_7s_4s_3s_4s_5s_4s_3s_1)^{-1}=-(s_1s_2s_3s_4s_5s_7)^{-1}.$$
This implies that the classes from the last line of the table are not real in $J_8$; thus there are 13 real classes in $J_8$.

Finally we examine the fifth column classifying the strongly real classes. Obviously involutions and the identity are strongly real. The only non-central involutions in $J_8$ correspond to 4-transpositions in $A_8$ (line 7 of the table). We can use this to rule out some cases: observe that
$$R_{A_8}((123))\cong (\langle (123)\rangle \times A_5):\langle (12)(45)\rangle,$$ hence $R_{A_8}((123))$ contains no 4-transpositions. We conclude that lines 3 and 9 do not correspond to strongly real classes. Similarly 
$$R_{A_8}((12345)) \cong (\langle (12345)\rangle \times \langle(678)\rangle):\langle (25)(34)\rangle,$$
and, again, this contains no 4-transpositions. Hence line 6 does not correspond to strongly real classes.

Now the group $R_{A_8}((123)(456))$ contains a 4-transposition, $(16)(25)(34)(78)$, which reverses $(123)(456)$. Since $(123)(456)$ lifts to elements of different orders, we conclude that they must be strongly real. In other words, line 4 corresponds to strongly real classes.

Next consider line 2 and take $g\in J_8$ which projects onto $gZ=(12)(34)$. Let $H$ be the group of even permutations of the set $\{5,6,7,8\}$ (so $H\cong A_4$). Then 
\begin{equation*}
 \begin{aligned}
C_{J_8}(g)/Z &= \{(1), (12)(34), (13)(24), (14)(23)\}\times H \\
C_{A_8}(gZ) &= C_{J_8}(g)/Z. \langle (12)(56)\rangle.
 \end{aligned}
\end{equation*}
Any element that reverses $g$ must centralize $gZ$. However all of the 4-transpositions that lie in $C_{A_8}(gZ)$ are contained in $C_{J_8}(g)/Z$; hence we conclude that $g$ is not strongly real.

We move on to line 5 and take $g\in J_8$ which projects onto $gZ=(1234)(56)$. Then 
\begin{equation*}
 \begin{aligned}
C_{J_8}(g)/Z &= \langle (1234)(56)\rangle, \\
C_{A_8}(gZ) &= \langle (1234)(56), (1234)(78) \rangle, \\
R_{A_8}(gZ) &= \langle (1234)(56), (1234)(78), (14)(23) \rangle.
 \end{aligned}
\end{equation*}
There are four cosets of $C_{J_8}(g)/Z$ in $R_{A_8}(gZ)$, two of which reverse $g$. Only one of these cosets contains a 4-transposition. Now observe that
$$(s_1s_5s_2s_1s_2s_3s_2s_1)(s_1s_2s_3s_5)(s_1s_5s_2s_1s_2s_3s_2s_1)^{-1} = -(s_1s_2s_3s_5)^{-1}.$$
Thus the coset containing a 4-transposition does not reverse $g$, and we conclude that line 5 does not correspond to a real class in $J_8$.

Finally consider line 8 and take $g\in J_8$ which projects onto $gZ=(1234)(5678)$ in $A_8$. Observe that
\begin{equation*}
 \begin{aligned}
C_{J_8}(g)/Z &= \langle (1234)(5678)\rangle, \\
C_{A_8}(gZ) &= \langle (1234)(5678), (1234)(8765)\rangle: \langle (15)(26)(37)(48)\rangle, \\
R_{A_8}(gZ) &= (C_{A_8}(gZ)):\langle (15)(26)(37)(48)\rangle.  
 \end{aligned}
\end{equation*}
Set $H=\langle (1234)(5678)\rangle;$ then $H$ has four cosets in $R_{A_8}(gZ)$, all of which contain 4-transpositions. One of these cosets must lift to the set of elements $\{h: hgh^{-1}=g^{-1}\}$, and we conclude that line 8 does correspond to a strongly real class in $J_8$.
\end{proof}

\subsection{Covers of $PSL_3(4)$}

Let $G=PSL_3(4)$, and let $H$ be a quasi-simple cover of $G$. Observe that $G$ contains a single conjugacy class of involutions; what is more \cite[Proposition 6.4.1]{gls} implies that an involution $g\in G$ lifts to an involution $h\in H$. 

Now, using our work above, we calculate that $G$ contains eight real classes (containing elements of order 1,2,3,4,4,4,5, and 5) and they are all strongly real. Furthermore $Z(H)\leq C_4 \times C_{12}$ and Lemma \ref{l: lifting} allows us to assume that $Z(H)$ is a non-trivial subgroup of $C_4\times C_4$. 

Thus there are seven covers of $L_3(4)$ which need to be addressed: 
\begin{equation*}
 \begin{aligned}
&PSL_3(4), \, 2.PSL_3(4), \, E_4.PSL_3(4), \, 4_1.PSL_3(4), \, 4_2.PSL_3(4), \\
&(E_4 4_1).PSL_3(4), \,(C_4\times C_4).PSL_3(4).  
 \end{aligned}
\end{equation*}

Note that by $E_4$ we mean an elementary-abelian group of order $4$; by $4_1$ and $4_2$ we mean quotients of $M(G)$ by cyclic groups of order $4$ that lie in $M(G)$ and are not in the same orbit of $Out(PSL_3(4))$. That this list of covers is comprehensive follows easily from \cite[Theorem 6.3.1]{gls} and \cite[Lemma 2.3(i), p.463]{gls6}.

Before we proceed with our analysis we need to establish some notation. Let $P$ be a Sylow $2$-subgroup of $PSL_3(4)$; in particular assume that
$$P = \left\{\left(\begin{matrix}
1 & a & b \\
0 & 1 & c \\
0 & 0 & 1 \end{matrix}\right) \mid a,b,c\in\mathbb{F}_4 \, \right\}.$$
Write $P_H$ for the Sylow $2$-subgroup of $H$ that projects onto $P$. For $h\in H$, define
$$Z_h=\{z\in Z \mid h_1hh_1^{-1} = hz \textrm{ for some } h_1\in P_H\}.$$
Observe that this is a subgroup of $Z$ and, that $|Z_h|=|C_P(g): C_{P_H}(h)/Z|$.

\begin{proposition}
Let $H$ be a quasi-simple cover of $G=PSL_3(4)$ with centre $Z$, a $2$-group. Suppose that $g=hZ\in G=H/Z$, with $g$ real in $G$ of order $d$.
\begin{enumerate}
\item If $d$ is odd, then $h$ is real if and only if the order of $h$ is $d$ or $2d$. What is more $h$ is strongly real if and only if $h$ is real.
\item If $d=2$, then $h$ is strongly real. 
\end{enumerate}
\end{proposition}

Note that we are not dealing with the case when $d=4$. We address this situation in the next proposition, using information from the Atlas \cite{atlas}. 

\begin{proof}
{\bf Suppose first that $d$ is odd}. Then the set $hZ$ generates a cyclic subgroup of $H$, and we may take the order of $h$ to equal $d$. If $d=1$ then $h$ is central and is (strongly) real if and only if $h^2=1$. 

Now suppose that $d>1$. Let $g_1\in G$ satisfy $g_1gg_1^{-1}=g^{-1}$ and $g_1=h_1Z$ for $h\in H$. Then $h_1hh_1^{-1}=h^{-1}z$ for some $z\in Z$. Since $h$ has odd order, we conclude that $z=1$. In other words $h_1hh_1^{-1}=h^{-1}$ and $h$ is real. More generally this implies that
$$h_1(hz)h_1^{-1}=h^{-1}z.$$
Thus $hz$ is real if and only if $z^2=1$. Furthermore, since $g$ is strongly real we may take $g_1$ to be an involution. By \cite[Proposition 6.4.1]{gls} we may therefore take $h_1$ to be an involution. Thus if $hz$ is real then $hz$ is strongly real; we have proved (1).

{\bf Suppose that $d=2$}. We take $h$ to have order $2$, which we may do by \cite[Proposition 6.4.1]{gls}.  We take
$$g=\left(\begin{matrix}
1 & 0 & x \\
0 & 1 & 0 \\
0 & 0 & 1 \end{matrix}\right)\in Z(P),$$
where $x\in\mathbb{F}_4^*$. 

We may assume that $H\cong (C_4\times C_4).PSL_3(4)$ and suppose that $Z(H)=\langle z_1, z_2\rangle$, so $z_1$ and $z_2$ are elements of order $4$. \cite[Lemma 2.3(c), p.463]{gls6} implies that $|C_P(g): C_{P_H}(h)/Z|=4$ and so $|Z_h|=4$. Since $h$ has order 2, this implies that $Z_h=\{1, z_1^2, z_2^2, z_1^2z_2^2\}$. The elements of $hZ$ can therefore be written in subsets of conjugate elements as follows:
\begin{equation*}
 \begin{aligned}
&\{h, hz_1^2, hz_2^2, hz_1^2z_2^2\}, \, \{hz_1, hz_1^3, hz_1z_2^2, hz_1^3z_2^2\}, \\
&\{hz_2, hz_1^2z_2, hz_2^3, hz_1^2z_2^3\}, \, \{hz_1z_2, hz_1^3z_2, hz_1z_2^3, hz_1^3z_2^3\}. 
 \end{aligned}
\end{equation*}
In particular these elements are all real. 

We need to establish that these elements are, in fact, strongly real. Observe first that if $h_1^2=hz$, with $h$ of order $2$ and $z$ central, then 
$$hh_1=h_1^{-1}z=h_1h.$$
In other words, an element $h_1$ satisfying $h_1^2=hz$ commutes with $h$. This, along with \cite[Lemma 2.3(c), p. 463]{gls6}, implies that $C_{P_H}(h)/Z$ is isomorphic to the group
\begin{equation*}
 \begin{aligned}
C &= \left\langle \left(\begin{matrix}
1 & a & b \\
0 & 1 & c \\
0 & 0 & 1 \end{matrix}\right) \mid a,b,c\in\mathbb{F}_4, ac=x \right\rangle \\
&= \left\{\left(\begin{matrix}
1 & a & b \\
0 & 1 & c \\
0 & 0 & 1 \end{matrix}\right) \mid a,b,c\in\mathbb{F}_4, ac=x \textrm{ or } a=c=0\, \right\}\cong C_2\times Q_8. 
 \end{aligned}
\end{equation*}

Now let $g_1$ be some element of $P$ such that $\langle C, g_1\rangle$ is a degree $2$ extension of $C_{P_H}(h)/Z$. Then 
$$g_1=\left(\begin{matrix}
1 & a & b \\
0 & 1 & c \\
0 & 0 & 1 \end{matrix}\right), \textrm{ for some } a,b,c\in\F_4,\textrm{ with } (a,c)\neq (0,0).$$
If $a$ or $c$ are equal to $0$ then $g_1^2=1$ and this extension is split. If $a\neq 0\neq c$, then observe that there exists $g_0\in C$ such that
$$g_0=\left(\begin{matrix}
1 & a & 0 \\
0 & 1 & a^{-1}x \\
0 & 0 & 1 \end{matrix}\right),$$
and $g_0g_1$ is an involution; thus, again, the extension is split. Thus any degree $2$ extension of $C_{P_H}(h)/Z$ in $P$ is split. 

Then, since $C_{P_H}(hz)/Z= C_{P_H}(h)/Z$ for any $z\in Z$, \cite[Proposition 6.4.1]{gls} implies that $R_P(hz)$ is a split extension of $C_P(hz)$. In other words, $hz$ is strongly real, as required.
\end{proof}

We must now examine those elements $h\in H$ for which $g=hZ$ is an element of order $4$ in $PSL_3(4)$. Note that there are three conjugacy classes of elements of order $4$ in $PSL_3(4)$. They are fused by an outer automorphism of $PSL_3(4)$ and intersect $P$ in sets which we label $C_k, k\in\mathbb{F}_4^*$:
$$C_k=\left\{\left(\begin{matrix}
1 & a & b \\
0 & 1 & c \\
0 & 0 & 1 \end{matrix}\right) \mid a,b,c\in\mathbb{F}_4, ac^{-1}=k\, \right\}.$$
If $Z(H)=E_4$ or $Z(H)=C_4\times C_4$ then the set of conjugacy classes in $H$ that project onto $C_k$ is mapped, via an outer automorphism of $H$, to the set of conjugacy classes in $H$ that project onto $C_{k'}$ for $k'\neq k$ \cite[Table 6.3.1]{gls}.

\begin{proposition}
Let $H$ be a quasi-simple cover of $G=PSL_3(4)$ with centre $Z$, a $2$-group. Suppose that $g=hZ\in G=H/Z$, with $g$ real in $G$ of order $4$. Suppose that $g$ lies in the set $C_k$ for some $k\in\F_4^*$.
\begin{enumerate}
\item If $Z(H)=C_2$ or $Z(H)=4_1$, then all of the elements in $hZ$ are real.
\item If $Z(H)=4_2$, then the number of real elements in $hZ$ depends on $k$. For two values of $k$, every element in $hZ$ is real; for the third, precisely half of the elements in $hZ$ are real.
\item If $Z(H)=C_4\times C_4$, then precisely half of the elements in $hZ$ are real.
\item If $Z(H)=E_44_1$, then the number of real elements in $hZ$ depends on $k$. For two values of $k$, every element in $hZ$ is real; for the third, precisely half of the elements in $hZ$ are real.
\item If $Z(H)=E_4$, then all of the elements in $hZ$ are real.
\item If $Z(H)$ is non-cyclic, then none of the elements in $hZ$ are strongly real.
\item If $Z(H)=C_4$, then the number of real elements in $hZ$ depends on $k$. For two values of $k$, none of the elements in $hZ$ are strongly real; for the third, precisely half of the elements in $hZ$ are strongly real.
\item If $Z(H)=C_2$, then the number of real elements in $hZ$ depends on $k$. For two values of $k$, none of the elements in $hZ$ are strongly real; for the third, all of the elements in $hZ$ are strongly real.
\end{enumerate}
\end{proposition}
\begin{proof}
Statements (1) and (2) follow immediately from \cite[p.28]{atlas}; we therefore start with (3). Throughout the proof we will refer to the {\it universal $2$-cover} as $H_U\cong (C_4\times C_4).PSL_3(4)$. To begin we take $H=H_U$, $x\in\mathbb{F}_4^*$, and set
\begin{equation*}
\begin{aligned}
g=\left(\begin{matrix}
1 & 1 & 0 \\
0 & 1 & x \\
0 & 0 & 1 \end{matrix}\right), \textrm{ thus } g^2=\left(\begin{matrix}
1 & 0 & x \\
0 & 1 & 0 \\
0 & 0 & 1 \end{matrix}\right).
\end{aligned}
\end{equation*}
This implies that
\begin{equation*}
C_P(g) = \left\{\left(\begin{matrix}
1 & a & b \\
0 & 1 & c \\
0 & 0 & 1 \end{matrix}\right) \mid a,b,c\in\mathbb{F}_4, c=ax\, \right\}\cong C_4\times C_4. 
\end{equation*}
Clearly $C_{P_H}(h)\geq \langle h\rangle Z(H)$ and so $C_{P_H}(h).Z$ has index at most $4$ in $C_P(g)$. Now we observe that one of the classes of order $4$ in $PSL_3(4)$ lifts to four separate conjugacy classes in $4_2.PSL_3(4)$ \cite[p.28]{atlas}. We conclude that $C_{P_H}(h)$ has index at least $4$ in $C_P(g)$; thus $C_{P_H}(h)/Z=\langle h\rangle Z/Z$, and $|C_P(g): C_{P_H}(h)/Z|=4=|Z_h|.$

Now write $Z=\langle z_1\rangle \times \langle z_2\rangle$. Suppose that $Z_h$ is elementary abelian; then
$$Z_h=\{1,z_1^2, z_2^2, z_1^2z_2^2\}.$$
Now there is an element in $C_{P_H}(h^2)/Z$ that conjugates $g$ to $g^{-1}$; hence it must map $h$ to an element of $h^{-1}Z_h$ (since $h^{-1}Z_h$ contains all elements of $h^{-1}Z$ whose square is equal to $h^2$). Since $Z_{h^{-1}}$ contains $h^{-1}$ we conclude that $h$ is real; indeed, all elements in $hZ$ are real. Since covers with cyclic centre are epimorphic images of covers with non-cyclic centre, the same conclusion will follow if $Z(H)$ is cyclic. This contradicts statements (1) and (2).

Hence we conclude that $Z_h$ is cyclic; write $Z_h=\{1, z_1, z_1^2, z_1^3\}$. Now there is an element $h_1\in C_{P_H}(h^2)/Z$ that conjugates $g$ to $g^{-1}$; hence we conclude that $Z_h=Z_{h^{-1}}$. We have several cases to consider:
\begin{enumerate}
\item If $h_1$ conjugates $hZ_h$ to $h^{-1}z_2^3Z_h$, then relabel so that $z_2^3$ becomes $z_2$; then we lie in the next case. 
\item If $h_1$ conjugates $hZ_h$ to $h^{-1}z_2Z_h$, then $hZ\cup h^{-1}Z$ splits into four sets of conjugate elements, with elements from distinct sets non-conjugate:
$$(hZ_h\cup h^{-1}z_2Z_h),\, z_2(hZ_h\cup h^{-1}z_2Z_h),\, z_2^2(hZ_h\cup h^{-1}z_2Z_h),\, z_2^3(hZ_h\cup h^{-1}z_2Z_h).$$
We conclude that none of these elements are real. Moreover, for $h_1\in hZ$, we find that $h_1\langle z_1\rangle$ is not real in $H/\langle z_1\rangle$. However $H/\langle z_1 \rangle$ is a cover of $PSL_3(4)$ with cyclic center. This contradicts statements (1) and (2).
\item If $h_1$ conjugates $hZ_h$ to $h^{-1}z_2^2Z_h$, then the set of conjugates of $h$ in $h^{-1}Z$ is equal to $\{h^{-1}z_2^2, h^{-1}z_1z_2^2, h^{-1}z_1^2z_2^2, h^{-1}z_1^3z_2^2\}$. In this case relabel so that $h$ becomes $hz_2$; then we lie in the next case.
\item The set of conjugates of $h$ in $h^{-1}Z$ is equal to $\{h^{-1}, h^{-1}z_1, h^{-1}z_1^2, h^{-1}z_1^3\}$. 
\end{enumerate}
Thus, provided we label appropriately, the following elements are all conjugate:
$$hZ_h\cup h^{-1}Z_{h} = \{h, hz_1, hz_1^2, hz_1^3, h^{-1}, h^{-1}z_1, h^{-1}z_1^2, h^{-1}z_1^3\}.$$
Similarly the following sets consist of conjugate elements:
$$z_2(hZ_h\cup h^{-1}Z_{h}), z_2^2(hZ_h\cup h^{-1}Z_{h}), z_2^3(hZ_h\cup h^{-1}Z_{h}).$$
Thus, of all elements in $hZ$, precisely the elements in the following sets are real in $H$:
$$hZ_h, z_2^2hZ_h.$$

We have proved (3). To prove (4) and (5) we examine the following sets of conjugate elements in $H_U$:
$$hZ_h\cup h^{-1}Z_{h}, z_2(hZ_h\cup h^{-1}Z_{h}), z_2^2(hZ_h\cup h^{-1}Z_{h}), z_2^3(hZ_h\cup h^{-1}Z_{h}).$$ 
Consider $H = H_U/ Z_1$ where $Z_1$ is a central subgroup of $H_U$. The following table lists those elements $h_1\in hZ$ for which $h_1Z_1$ is real in $H$:

\begin{center}
\begin{tabular}{|c|c|c|}
\hline
$Z_1$ & $H$ & Real elements\\
\hline
$\langle z_1^2\rangle$ & $(E_44_1).PSL_3(4)$ & $Z_h, \, z_2^2 Z_h$ \\
$\langle z_2^2\rangle$ & $(E_44_1).PSL_3(4)$ & $Z_h, \, z_2Z_h, \, z_2^2Z_h, z_2^3Z_h$ \\
$\langle z_1^2z_2^2\rangle$ & $(E_44_1).PSL_3(4)$ & $Z_h, \, z_2Z_h, z_2^2Z_h, \, z_2^3Z_h$ \\
$\langle z_1^2, z_2^2\rangle$ & $E_4.PSL_3(4)$ & $Z_h, \, z_2Z_h, \, z_2^2Z_h, \, z_2^3Z_h$ \\
\hline
\end{tabular}
\end{center}

This yields (4) and (5). Note that we have three entries corresponding to $H=(E_44_1).PSL_3(4)$ as the order $3$ automorphisms of $PSL_3(4)$ do not lift to this group.

To prove the remaining statements we must determine when $h$ is strongly real.  Suppose first that $Z$ is cyclic and non-trivial. We start by considering the $4$-covers of $PSL_3(4)$; let $Z=\langle z \rangle$ and let $Y$ be the pre-image of $Z(P)$ in $P_H$; then \cite[Lemma 2.2, p.463]{gls6} implies that $Y=X\times Z$ where $X\cong C_2\times C_2$. Furthermore \cite[Lemma 2.3(e), p.463]{gls6} implies that $C_{P_H}(Y)/Z\cong C_4\times C_4$ and so contains an element $h$ such that $hZ$ is an element of order $4$ in $PSL_3(4)$. In particular $C_{P_H}(h)>Y$.

Now suppose that $H\cong 4_1.PSL_3(4)$; then $C_{P_H}(h)/Z$ is a proper subgroup of $C_P(g)$ (otherwise, $hZ$ intersects four distinct conjugacy classes of $4_1.PSL_3(4)$, and $Z_h=\{1\}$; this is impossible \cite[p. 24]{atlas}). Thus $C_{P_h}(h)=\langle h \rangle Y$, which has index $2$ in $C_P(g)$. This implies, firstly, that $hZ$ intersects two conjugacy classes, call them $C_1$ and $C_2$, in $H$. It implies, secondly, that $Z_h=\{1,z^2\}=Z_{h^{-1}}$ and, since $h$ is real, the elements $h, h^{-1}, hz^2$ and $h^{-1}z^2$ are all conjugate.

Now the set of elements in $P$ that reverse $g$ is equal to
$$R=\left\{\left(\begin{matrix}
1 & a & b \\
0 & 1 & (a+1)x \\
0 & 0 & 1 \end{matrix}\right) \mid a,b\in\mathbb{F}_4 \, \right\}.$$
In addition observe that $C= C_{P_H}(h)/Z\ = \langle g, Z(P) \rangle$. Consider degree 2 extensions of $C$ of the form $\langle C, r\rangle$ for some $r\in R$. There are two such extensions, one split (when the element $r$ has $a=1$ or $a=0$ in the matrix form given above) and the other non-split.

Thus $h$ is mapped by an involution to precisely one of either $h^{-1}$ or $h^{-1}z^2$. This implies that $h$ is strongly real if and only if $hz$ is not strongly real. Thus we conclude that precisely two of the elements in $hZ$ are strongly real in $4_1.PSL_3(4)$. In particular not all real elements are strongly real in $4_1.PSL_3(4)$.

We return to the situation where $H=H_U$ and write $Z(H)=\langle z_1, z_2 \rangle$. As before we choose $h$ so that 
$$Z_h=\{1, z_1, z_1^2, z_1^3\}.$$
Again the following sets consist of conjugate elements:
$$hZ_h\cup hZ_{h^{-1}}, z_2(hZ_h\cup hZ_{h^{-1}}), z_2^2(hZ_h\cup hZ_{h^{-1}}), z_2^3(hZ_h\cup hZ_{h^{-1}}).$$
Thus, in $H/\langle z_1 z_2^2 \rangle \cong 4_1.PSL_3(4)$, the set $hZ$ splits into two conjugacy classes; then \cite[p. 24]{atlas} implies that these conjugacy classes must coincide with $C_1$ and $C_2$ described above.

Now define groups $R<S< G$ as follows:
\begin{equation*}
\begin{aligned}
R &= \left\{\left(\begin{matrix}
1 & a & b \\
0 & 1 & c \\
0 & 0 & 1 \end{matrix}\right) \mid a,b,c\in\mathbb{F}_4, a\in\{0,1\}\, c\in \{0,x\} \right\}\cong C_4\times C_4; \\
S &= \left\{\left(\begin{matrix}
1 & a & b \\
0 & 1 & c \\
0 & 0 & 1 \end{matrix}\right) \mid a,b,c\in\mathbb{F}_4, (a,c)\in\{(0,0), (1,x)\} \right\}\cong C_4\times C_2. 
\end{aligned}
\end{equation*}
Let $R_H$ (resp. $S_H$) be the pre-image of $R$ (resp. $S$) in $H$; then $R_H$ is a degree $4$ extension of $C_{P_H}(h)$. In addition, firstly, all of the involutions that reverse $g$ are contained in $R$. Secondly, $S_H$ is a subgroup of $R_H$ of index $2$, and $S_H/Z$ centralizes $g$. Thus
$$\{z\in Z \mid h_1hh_1^{-1}=hz \textrm{ for some } h_1\in S_H\},$$
is a subgroup of $Z_h$ size $2$; it must equal $\{1, z_1^2\}$. 

The set of $R$-conjugates of $h$ can take on two possible shapes; the first possibility is that the $R$-conjugates of $h$ are equal to 
$$h^R=\{h, hz_1^2, h^{-1}, h^{-1}z_1^2\}.$$
The strongly real elements of $hZ$ are, then, all elements $hz$ such that the set $zh^R$ contains $(hz)^{-1}$. A quick calculation demonstrates that the elements satisfying this requirement are precisely the real elements in $hZ$. Thus all real elements in $H$ are strongly real. Clearly the same result applies to all epimorphic images of $H$, which contradicts our earlier calculations in $4_1.PSL_3(4)$.

The second possibility is that the $R$-conjugates of $h$ are equal to $\{h, hz_1^2, h^{-1}z_1, h^{-1}z_1^3\}$. In this case half of the elements in $hZ$ are strongly real in $H/\langle z_1 z_2^2 \rangle$, which is consistent with our calculations above. It immediately follows that none of the elements in $hZ$ are strongly real in $H_U$.

To complete our analysis we consider, as before, $H = H_U/ Z_1$ where $Z_1$ is a central subgroup of $H_U$. The following table lists those choices of $Z_1$ for which $hZ$ contains any strongly real elements; in each case the table lists those elements $h_1\in hZ$ for which $h_1Z_1$ is strongly real in $H$; we write $Y_h$ for the set $\{1, z_1^2\}$:

\begin{center}
\begin{tabular}{|c|c|c|}
\hline
$Z_1$ & $H$ & Strongly real elements\\
\hline
$\langle z_1z_2^2\rangle$ & $4_1.PSL_3(4)$ & $hz_2Y_h, \, hz_1z_2Y_h, \, hz_2^3Y_h, \, hz_1z_2^3Y_h$ \\
$\langle z_1\rangle$ & $4_2.PSL_3(4)$ & $hY_h, \, hz_1Y_h, \, hz_2Y_h, \, hz_1z_2^2Y_h$ \\
$\langle z_1, z_2^2\rangle$ & $2.PSL_3(4)$ & $hZ$ \\
\hline
\end{tabular}
\end{center}

Statements (6), (7) and (8) follow immediately from the table.
\end{proof}

\section{Some small rank calculations}\label{s: small}

We list our results for $n\leq 6$ (excluding the exceptional cases covered in the previous section). For $k$ a positive integer we write $\delta_k$ to mean $(q-1,k)$.

\subsection{$q$ even}

In this case the number of real and strongly real elements always coincides. Furthermore the counts for $GL_n(q)$ and $PGL_n(q)$ coincide, as do the counts for $SL_n(q)$ and $PSL_n(q)$. Thus for each $n$ we simply record the number of real classes in $GL_n(q)$ and in $SL_n(q)$:

\begin{center}
\begin{tabular}{|c|c|c|}
\hline
$n$ & $GL_n(q)$ & $SL_n(q)$ \\
\hline
2 & $q+1$ & $q+1$ \\
3 & $q+2$ & $q+1+\delta_3$ \\
4 & $q^2+2q+2$ & $q^2 + 2q+2$ \\
5 & $q^2+3q+3$ & $q^2+3q+2+\delta_5$ \\
6 & $q^3+2q^2+4q+4$ & $q^3+2q^2+ (3+\delta_3)q + 3+\delta_3$ \\
\hline
\end{tabular}
\end{center}

\subsection{$q$ odd, $n\not\equiv 2\pmod 4$}

In this case the number of real and strongly real elements always coincides. Thus we record only counts for real classes.

\begin{center}
\begin{tabular}{|c|c|c|c|c|}
\hline
$n$ & $GL_n(q)$ & $SL_n(q)$ & $PGL_n(q)$ & $PSL_n(q)$ \\
\hline
3 & $2q+6$ & $q+2+\delta_3$ & $q+3$ & $q+2+\delta_3$\\
4 & $q^2+4q+9$ & $q^2+4q+4+2\delta_4$ & $q^2+3q+5$ & $\ast$\\
5 & $2q^2+8q+14$ & $q^2+4q+6+\delta_5$ & $q^2+4q+7$ & $q^2+4q+6+\delta_5$ \\
\hline
\end{tabular}
\end{center}

We must fill in the values for $PSL_4(q)$. The number of real classes are as follows:
$$\left\{ \begin{array}{ll}
\frac12q^2+\frac52q+3+\delta_4, & q\equiv 1 \pmod 4  \\
q^2+3q+3+\delta_4, & q\equiv 3 \pmod 4  \\
\end{array}\right.
$$

\subsection{$q$ odd, $n=2$}

All real classes are strongly real in both $GL_2(q)$ and $PGL_2(q)$. There are $q+3$ such classes in $GL_2(q)$ and $q+2$ such classes in $PGL_2(q)$.

In $SL_2(q)$ there are precisely $2$ strongly real classes. The number of real classes is as follows:
$$\left\{ \begin{array}{ll}
q+4, & q\equiv 1 \pmod 4  \\
q, &   q\equiv 3 \pmod 4. \\
\end{array}\right.
$$

In $PSL_2(q)$ all real classes are strongly real. The number of real classes is as follows:
$$\left\{ \begin{array}{ll}
\frac12 q+\frac 52, & q\equiv 1 \pmod 4  \\
\frac12q+\frac12, &   q\equiv 3 \pmod 4. \\
\end{array}\right.
$$

\subsection{$q$ odd, $n=6$}

All real classes are strongly real in both $GL_6(q)$ and $PGL_6(q)$. There are $q^3+4q^2+13q+22$ such classes in $GL_6(q)$ and $q^3+3q^2+9q+12$ such classes in $PGL_6(q)$.

In $SL_6(q)$ there are precisely $4q^2+8q+12+2\delta_3$ strongly real classes. The number of real classes is as follows:
$$\left\{ \begin{array}{ll}
q^3+3q^2+(9+\delta_3)q+14+4\delta_3, & q\equiv 1 \pmod 4  \\
q^3+3q^2+(5+\delta_3)q+6, &   q\equiv 3 \pmod 4. \\
\end{array}\right.
$$

In $PSL_6(q)$ there are $q^3+2q^2+(7+\delta_3)q+7+2\delta_3$ real classes when $q\equiv 1\pmod 4$; what is more all of these classes are strongly real. When $q\equiv 3\pmod 4$ there are $\frac12q^3+2q^2+(3+\frac12\delta_3)q+\frac72+\frac12\delta_3$ real classes and all real classes are strongly real except for some of type $\nu=1^4 2^1$. Thus the total number of strongly real classes is $\frac12q^3+2q^2+(3+\frac12\delta_3)q+\frac72+\frac12\delta_3-x$ where $x=\frac12(q^2-q)$ is the number of monic self-reciprocal polynomials of degree $4$ with constant term $1$ which are either irreducible or else are the product of two non-self-reciprocal degree $2$ irreducible polynomials.

\section{Further work}\label{s: further}

It is very natural to ask if the real and strongly real classes can be counted in other families of finite groups of Lie type. In fact the work of Macdonald extends naturally (as he explains in \cite{macdonald}) to the unitary groups, so this is the natural next step.

By counting real conjugacy classes for a finite group $G$ we are also, of course, counting real irreducible representations for $G$. The question is: can we now construct these representations? For the case of $GL_n(q)$ we hope to use Green's classical method to do just that \cite{green}. For other finite groups of Lie type this is likely to be very difficult, and requires an understanding of the Deligne-Lusztig theory.

Note finally that real irreducible characters come from two different kinds of irreducible representation, the orthogonal and the symplectic. It is not clear whether there is any such division for the number of real conjugacy classes, however we can make some interesting observations. For instance, note that for $GL_n(q)$, and for  $SL_n(q)$ with $n\not\equiv 2 \imod 4$, all real conjugacy classes are strongly real; it turns out that in these cases the self dual representations are orthogonal, i.e., the real characters actually come from orthogonal representations (c.f. \cite{gow2, ohmori} where it is shown that the Schur index for the characters of above mentioned groups are $1$). 

This correspondence between strongly real classes and orthogonal representations has been studied from a variety of angles -- see Gow's work on 2-regular structure  \cite{gow3}, and Prasad's work on groups of Lie type and $p$-adic groups \cite{prasad1, prasad2}. Nonetheless, although the correspondence can be seen to hold in particular cases (and not in others), it is unclear how general a phenomenon it really is.

\bibliographystyle{amsalpha}
\bibliography{paper}

\newcommand{\etalchar}[1]{$^{#1}$}
\providecommand{\bysame}{\leavevmode\hbox to3em{\hrulefill}\thinspace}
\providecommand{\MR}{\relax\ifhmode\unskip\space\fi MR }
\providecommand{\MRhref}[2]{%
  \href{http://www.ams.org/mathscinet-getitem?mr=#1}{#2}
}
\providecommand{\href}[2]{#2}
\begin{thebibliography}{CCN{\etalchar{+}}85}

\bibitem[Asc00]{aschbacher3}
Michael Aschbacher, \emph{Finite group theory}, Cambridge studies in advanced
  mathematics, no.~10, Cambridge University Press, 2000.

\bibitem[CCN{\etalchar{+}}85]{atlas}
J.H. Conway, R.T. Curtis, S.P. Norton, R.A. Parker, and R.A. Wilson,
  \emph{Atlas of finite groups}, Oxford University Press, 1985.

\bibitem[GLS94]{gls}
Daniel Gorenstein, Richard Lyons, and Ronald Solomon, \emph{The classification
  of the finite simple groups, number 3}, Mathematical Surveys and Monographs,
  vol.~40, American Mathematical Society, 1994.

\bibitem[GLS05]{gls6}
\bysame, \emph{The classification of the finite simple groups. {N}umber 6.
  {P}art {IV}}, Mathematical Surveys and Monographs, vol.~40, American
  Mathematical Society, Providence, RI, 2005, The special odd case.

\bibitem[Gow76a]{gow2}
R.~Gow, \emph{Schur indices of some groups of {L}ie type}, J. Algebra
  \textbf{42} (1976), no.~1, 102--120.

\bibitem[Gow76b]{gow3}
Roderick Gow, \emph{Real-valued characters and the {S}chur index}, J. Algebra
  \textbf{40} (1976), no.~1, 258--270.

\bibitem[Gow81]{gow}
R.~Gow, \emph{The number of equivalence classes of nondegenerate bilinear and
  sesquilinear forms over a finite field}, Linear Algebra Appl. \textbf{41}
  (1981), 175--181.

\bibitem[Gre55]{green}
J.~A. Green, \emph{The characters of the finite general linear groups}, Trans.
  Amer. Math. Soc. \textbf{80} (1955), 402--447.

\bibitem[Isa94]{isaacs}
I.~Martin Isaacs, \emph{Character theory of finite groups}, Dover Publications
  Inc., 1994.

\bibitem[Jac53]{jacobson2}
Nathan Jacobson, \emph{Lectures in abstract algebra. {V}ol. {II}. {L}inear
  algebra}, D. Van Nostrand Co., Inc., Toronto-New York-London, 1953.

\bibitem[KL90]{kl}
P.~Kleidman and M.~Liebeck, \emph{The subgroup structure of the finite simple
  groups}, London Mathematical Society Lecture Note Series, vol. 129, Cambridge
  University Press, Cambridge, 1990.

\bibitem[Lan02]{lang}
Serge Lang, \emph{Algebra}, third ed., Graduate Texts in Mathematics, vol. 211,
  Springer-Verlag, New York, 2002.

\bibitem[Leh75]{lehrer}
G.~I. Lehrer, \emph{Characters, classes, and duality in isogenous groups}, J.
  Algebra \textbf{36} (1975), no.~2, 278--286.

\bibitem[Mac81]{macdonald}
I.~G. Macdonald, \emph{Numbers of conjugacy classes in some finite classical
  groups}, Bull. Austral. Math. Soc. \textbf{23} (1981), no.~1, 23--48.

\bibitem[Ohm77]{ohmori}
Zyozyu Ohmori, \emph{On the {S}chur indices of {${\rm GL}(n, q)$} and {${\rm
  SL}(2n+1, q)$}}, J. Math. Soc. Japan \textbf{29} (1977), no.~4, 693--707.

\bibitem[Pra98]{prasad1}
Dipendra Prasad, \emph{On the self-dual representations of finite groups of
  {L}ie type}, J. Algebra \textbf{210} (1998), no.~1, 298--310.

\bibitem[Pra99]{prasad2}
\bysame, \emph{On the self-dual representations of a {$p$}-adic group},
  Internat. Math. Res. Notices (1999), no.~8, 443--452.

\bibitem[Spr94]{springer2}
T.A. Springer, \emph{Linear algebraic groups}, 1--123, Algebraic Geometry IV,
  Encyclopedia of Mathematical Sciences 55.

\bibitem[ST05]{singhthakur2}
Anupam Singh and Maneesh Thakur, \emph{Reality properties of conjugacy classes
  in algebraic groups}, Israel J. Math. \textbf{145} (2005), 157--192.

\bibitem[TZ05]{tz}
Pham~Huu Tiep and A.~E. Zalesski, \emph{Real conjugacy classes in algebraic
  groups and finite groups of {L}ie type}, J. Group Theory \textbf{8} (2005),
  no.~3, 291--315.

\bibitem[Vin04]{vinroot}
C.~Ryan Vinroot, \emph{A factorization in {${\rm GSp}(V)$}}, Linear Multilinear
  Algebra \textbf{52} (2004), no.~6, 385--403.

\bibitem[Wal80]{wall}
G.~E. Wall, \emph{Conjugacy classes in projective and special linear groups},
  Bull. Austral. Math. Soc. \textbf{22 no. 3} (1980), 339--364.

\bibitem[Won66]{wonen}
Mar{\'{\i}}a~J. Wonenburger, \emph{Transformations which are products of two
  involutions}, J. Math. Mech. \textbf{16} (1966), 327--338.

\end{thebibliography}

\end{document}